\documentclass[12pt,a4paper,reqno]{amsart}
\usepackage{amssymb}
\usepackage{latexsym}
\usepackage{exscale}
\usepackage{mathrsfs}

\newcommand{\R}{\mathbb R}

\newcommand{\B}{\mathcal B}

\newcommand {\SB} {{\mathbb B}}

\newcommand {\SR} {{\mathbb R}}

\numberwithin{equation}{section}
\newtheorem{theorem}{Theorem}[section]
\newtheorem{lemma}[theorem]{Lemma}

\newtheorem{corollary}[theorem]{Corollary}
\newtheorem{Remark}[theorem]{Remark}
\newtheorem{proposition}[theorem]{Proposition}
\newtheorem{definition}[theorem]{Definition}
\newtheorem{example}[theorem]{Example}

\begin{document}
\allowdisplaybreaks
\title{Best Approximation with wavelets in Weighted  Orlicz Spaces}
\author{Maria de Natividade}
\address{Maria de Natividade
\\
Departamento de Matem\'aticas
\\
Universidad Aut\'onoma de Ma\-drid
\\
28049 Madrid, Spain} \email{maria.denatividade@uam.es}
\begin{abstract}
Democracy functions of wavelet admissible bases are computed for weighted Orlicz Spaces $L^{\Phi}(w)$ in terms of the fundamental function of $L^{\Phi}(w).$ In particular, we prove that these bases are greedy in $L^{\Phi}(w)$ if and only if $L^{\Phi}(w) = L^{p}(w),\ 1<p<\infty.$ Also, sharp embeddings for the approximation spaces are given in terms of weighted discrete Lorentz spaces. For $L^{p}(w)$ the approximation spaces are identified with weighted Besov spaces.
\end{abstract}
\thanks{\textit{Acknowledgements} and \textit{Notes}: This research  was supported by Instituto Nacional de Bolsas de Estudo de Angola, INABE and by Grant MTM2007-60952 of Spain. The author wishes to
 thank J. M. Martell  and G. Garrig\'os for useful discussions on this topic.}
\date{\today}
 \maketitle
\section{Introduction}\label{secIntr}
 Let ($\mathbb{B},\|\cdot\|_{\mathbb{B}}$) be a
quasi-Banach space with a countable unconditional basis
$\mathcal{B} = \{e_{j}: j \in \mathbb{N}\}$; that is, every $x \in
\mathbb{B}$ can be uniquely represented as an unconditionally
convergent series $x = \sum_{j\in\mathbb{N}}s_{j}e_{j}$, for some
sequences of scalars $\{s_{j}: j\in\mathbb{N}\}$. Let $\Sigma_{N}$
denote the set of all elements $y \in \mathbb{B}$ with at most $N$
non-null coefficients in the basis representation $y =
\sum_{j\in\mathbb{N}}s_{j}e_{j}$. For $x \in \mathbb{B}$, the
$N$-{\bf term error of approximation} (with respect to
$\mathcal{B}$) is defined by
\begin{eqnarray}\label{errorApproximation}\sigma_{N}(x)_{\mathbb{B}}\equiv \inf_{y
\in\Sigma_{N}}\|x-y\|_{\mathbb{B}}.\end{eqnarray} Two main questions
in approximation theory concern the construction of efficient
algorithms for $N$-term approximation and the characterization of
the {\bf approximation spaces}
$\mathcal{A}_{q}^{\alpha}(\mathcal{B},\mathbb{B})$, which consists of all $x
\in\mathbb{B}$ such that the quantity
\begin{eqnarray}\label{ApproximationSpaceNorm}\|x\|_{\mathcal{A}_{q}^{\alpha}(\mathcal{B},\mathbb{B})}=\left\{%
\begin{array}{ll}
   \Big(\sum_{N\geq 1}(N^{\alpha}\sigma_{N}(x)_{\mathbb{B}})^{q}\frac{1}{N}\Big)^{\frac{1}{q}} , & \hbox{ if $0<q<\infty$;} \\
   \\
  \sup_{N\geq1}[N^{\alpha}\sigma_{N}(x)_{\mathbb{B}}]  , & \hbox{if $q = \infty$,} \\
\end{array}%
\right.\end{eqnarray} is finite. A computational efficient method to
produce $N$-term approximations, which has been widely investigate
in recent years, is the so called
\textbf{greedy algorithm} (see e.g \cite{Konyagin}). If $ x = \sum_{j\in\mathbb{N}}s_{j}e_{j}$
and we order the basis elements in such a way that
\begin{eqnarray}\label{nonincreaRearrang}\|s_{j_{1}}e_{j_{1}}\|_{\mathbb{B}}\geq\|s_{j_{2}}e_{j_{2}}\|_{\mathbb{B}}\geq\ldots\end{eqnarray}
(handling ties arbitrarily), the {\bf greedy algorithm of step} $N$
is defined by the correspondence
\begin{eqnarray}\label{Greedyalgorithm} x
=\sum_{j\in\mathbb{N}}s_{j}e_{j}\in\mathbb{B}\longrightarrow
G_{N}(x) =
\sum_{k=1}^{N}s_{j_{k}}e_{j_{k}}\in\Sigma_{N}.\end{eqnarray}
S.V. Konyagin and V. N. Temlyakov (\cite{Konyagin}) defined the basis
$\mathcal{B}$ to be \textbf{greedy} in   $(\mathbb{B},
\|\cdot\|_{\mathbb{B}})$  if the greedy algorithm is optimal in the
sense  that $G_{N}(x)$ is essentially the best $N$-term
approximation to $x$ using the basis vectors, i.e, there exists a
constant $C$ such that for all $x \in \mathbb{B}$ we have
\begin{eqnarray*}\|x - G_{N}(x)\|_{\mathbb{B}}\leq
C\sigma_{N}(x)_{\mathbb{B}},\quad N = 1, 2,\ldots.\end{eqnarray*}
Thus, for such bases the greedy algorithm produces an almost optimal
$N$-term approximation, which leads often to a precise
identification of the approximation spaces
$\mathcal{A}_{q}^{\alpha}(\mathcal{B},\mathbb{B})$. In \cite{Konyagin} greedy
basis in a  quasi-Banach space
($\mathbb{B},\|\cdot\|_{\mathbb{B}}$) are characterized  as those which
are unconditional and {\bf democratic}, the latter meaning that
there exists some constant $\Delta > 0$ such that
\begin{eqnarray*}\Big\|\sum_{j\in\Gamma}\frac{e_{j}}
{\|e_{j}\|_{\mathbb{B}}}\Big\|_{\mathbb{B}}\leq\Delta\Big\|
\sum_{j\in\Gamma'}\frac{e_{j}}{\|e_{j}\|_{\mathbb{B}}}\Big\|_{\mathbb{B}},
\end{eqnarray*}
holds for all finite sets of indices $\Gamma,\Gamma'\subset
\mathbb{N}$ with the same cardinality. Wavelet systems are well
known examples of greedy bases for many function   and distribution
spaces. Indeed, V.N. Temlyakov showed in \cite{T1} that the Haar basis
(and any wavelet system $L^{p}$-equivalent to it) is greedy in the
Lebesgue space $L^{p}([0,1])$ for $1<p<\infty$. When
wavelet have sufficient smoothness and decay, they are also greedy
bases for the more general Sobolev and Triebel-Lizorkin classes (see
e.g \cite{HJLY,GustavoSharpJackson}).

 The purpose of this paper is
to study the efficiency of wavelet greedy algorithms in the weighted  Orlicz spaces $L^{\Phi}(w)$ defined for functions on $\mathbb{R}^{d}$. In Theorem \ref{thmcaraOrlPes} (see section \ref{secPrelim} ) we show that wavelet bases are unconditional in weighted Orlicz spaces $L^{\Phi}(w)$ with nontrivial Boyd indices  for all $w\in A_{p^{\Phi}}(\mathbb{R}^{d}).$  We give in section \ref{SecDemocFunOrl} a simple proof of the fact  that admissible wavelet bases (see definition below) are not democratic in weighted Orlicz spaces $L^{\Phi}(w)$  if $L^{\Phi}(w) \neq L^{p}(w).$ \\

   In view of this result it have interest to ask how far wavelet bases  are from being democratic in $L^{\Phi}(w)\neq L^{p}(w).$  To quantify democracy of a basis $\mathcal{B}=\{e_{j}\}_{j\in\mathbb{N}}$ we shall study the following functions: $$h_{r}(N;\mathcal{B},\mathbb{B}) = \sup_{Card(\Gamma)=N}\Big\|\sum_{\gamma\in\Gamma}\frac{e_{\gamma}}{\|e_{\gamma}\|_{\mathbb{B}}}\Big\|_{\mathbb{B}}\quad and \quad h_{l}(N;\mathcal{B},\mathbb{B})= \inf_{Card(\Gamma)=N}\Big\|\sum_{\gamma\in\Gamma}\frac{e_{\gamma}}{\|e_{\gamma}\|_{\mathbb{B}}}\Big\|_{\mathbb{B}}$$ which we call right and left democracy functions of $\mathcal{B}$ (see also \cite{DKKT, KT, GHM}). Observe that a basis is democratic if and only if these two quantities are comparable for all $N\geq 1.$ Our main result gives a precise value (except for multiplicative constants) of these functions in terms of intrinsic properties of the space $L^{\Phi}(w).$ Namely, let $h_{\varphi}^{+}(t) = \sup_{s>0}\frac{\varphi(st)}{\varphi(s)}$ denote the dilation function associated with the \textit{fundamental function} $\varphi$ of $L^{\Phi}(w),$ and let $h_{\varphi}^{-}(t)$  be the same quantity with ``$\sup$'' replaced by ``$\inf$'' (see Section \ref{secPrelim} for precise definitions).  \begin{theorem}\label{thmPrinc}Let $L^{\Phi}(w)$ be a weighted Orlicz space with non trivial Boyd indices, $w\in A_{p^{\Phi}},$ a weight on $\mathbb{R}^{d},$ and $\mathcal{B} = \{\psi_{Q}: Q\in \mathcal{D}\}$ be an admissible wavelet basis. Then for all $\Gamma\subset\mathcal{D}$ \begin{eqnarray}\label{thmPrinc1}h_{r}(N; \mathcal{B},L^{\Phi}(w))\approx h_{\varphi}^{-}(N),\quad \quad h_{l}(N;\mathcal{B},L^{\Phi}(w))\approx h_{\varphi}^{+}(N).\end{eqnarray} (Here $p^{\Phi} = \frac{1}{I_{\varphi}},$ where $I_{\varphi}$ is the upper Boyd index of $L^{\Phi}(w).$ See definition of Boyd indices in subsection 2.1.)\end{theorem}  This result will have  applications in the study of  approximation spaces (defined using admissible wavelet basis) in weighted Orlicz spaces. We  take up  this task  in the section \ref{secIncluAppr}, where we investigate Jackson and Bernstein type estimates and corresponding inclusions for $N$-term approximation spaces. In the  $L^{p}$ case, these estimates are naturally given in terms of the class of discrete Lorentz spaces $\ell^{\tau,q}$ (see \cite{DPetrova,GustavoSharpJackson,GNielsen,HJLY,KerkiacharyanNon}). In the case of weighted Orlicz spaces we shall need \textit{weighted Lorentz sequence spaces} $\Lambda_{\eta}^{q}$, defined by   \begin{eqnarray}\label{defnLorDisSp}\Lambda_{\eta}^{q} = \Big\{\textbf{s}:\|\textbf{s}\|_{\Lambda_{\eta}^{q}} = \Big[\sum_{k\geq 1}(\eta_{k}|s_{k}^{\ast}|)^{q}\frac{1}{k}\Big]^{\frac{1}{q}}<\infty\Big\}.\end{eqnarray} where $\{s_{k}^{\ast}\}$ is the non-increasing rearrangement of $\textbf{s}$ and the weight $\eta = \{\eta_{k}\}$ is a fixed increasing and doubling sequence (see \cite{GHM}). In particular, $\Lambda_{\eta}^{q} = \ell^{\tau,q}$ when $\eta_{k} = k^{1/\tau}.$

   For $f\in L^{\Phi}(w),$ and $\mathcal{B} = \{\psi_{Q}: Q\in\mathcal{D}\}$ a wavelet basis in $L^{\Phi}(w),$   write $f=\sum_{Q\in\mathcal{D}}\langle f,\psi_{Q}\rangle\psi_{Q}.$ Then we define  $\Lambda_{\eta}^{q}(\mathcal{B},L^{\Phi}(w))$ as the set of all $f\in L^{\Phi}(w)$ such that the sequence $\{\|\langle f,\psi_{Q_{k}}\rangle\psi_{Q_{k}}\|_{L^{\Phi}(w)}: k\geq 1\}\in \Lambda_{\eta}^{q}$ and $$\|f\|_{\Lambda_{\eta}^{q}(\mathcal{B},L^{\Phi}(w))} = \Big\|\|\langle f,\psi_{Q_{k}}\rangle\psi_{Q_{k}}\|_{L^{\Phi}(w)}\Big\|_{\Lambda_{\eta}^{q}}$$ where $\|\langle f,\psi_{Q_{1}}\rangle\psi_{Q_{1}}\|_{L^{\Phi}(w)}\geq \|\langle f,\psi_{Q_{2}}\rangle\psi_{Q_{2}}\|_{L^{\Phi}(w)}\geq\ldots$ (handling ties arbitrarily).
   \begin{theorem}\label{thmInclAprOrlpesoSpa}Let $L^{\Phi}(w)$ be a weighted Orlicz space with Boyd indices $0<i_{\varphi}\leq I_{\varphi}<1,$ and $w\in A_{p^{\Phi}}$ a weight on $\mathbb{R}^{d},$  where $p^{\Phi} =\frac{1}{I_{\varphi}}.$ Then \begin{eqnarray}\label{thminclAppOrlPesoSpa1}\Lambda_{k^{\alpha}h_{r}(k)}^{q}(\mathcal{B},L^{\Phi}(w))\hookrightarrow \mathcal{A}_{q}^{\alpha}(\mathcal{B},L^{\Phi}(w))\hookrightarrow \Lambda_{k^{\alpha}h_{l}(k)}^{q}(\mathcal{B},L^{\Phi}(w)),\end{eqnarray} \end{theorem}These embeddings are optimal, in the sense that the largest and smallest weighted Lorentz spaces $\Lambda_{k^{\alpha}\eta(k)}^{q}(\mathcal{B},L^{\Phi}(w))$ that one can place on the left- and right-hand side of (\ref{thminclAppOrlPesoSpa1}) are respectively $\Lambda_{k^{\alpha}h_{r}(k)}^{q}(\mathcal{B},L^{\Phi}(w))$ and $\Lambda_{k^{\alpha}h_{l}(k)}^{q}(\mathcal{B},L^{\Phi}(w))$ (see section 4).  We point out that a  sufficient condition for these two spaces to be equal is that $h_{r}(N)\approx h_{l}(N),$ in which case the basis is necessarily democratic and $L^{\Phi}(w) = L^{p}(w)$ (see Lemma 5.2 in \cite{GHM}). Then Theorem \ref{thmInclAprOrlpesoSpa} leads the following identification of Approximation spaces for $L^{p}(w)$ in terms of classical Lorentz spaces. \begin{corollary}\label{coro-CarLpesAppr} Let $\alpha>0,  \ 1<p<\infty,\ 0<q\leq\infty,$  and  $w\in A_{p}$ a weight on $\mathbb{R}^{d}.$ Then, for a wavelet basis $\mathcal{B},$ we have \begin{eqnarray}\label{coro-CarLpesAppr1}\mathcal{A}_{q}^{\alpha}(\mathcal{B},L^{p}(w)) = \ell^{\tau,q}(\mathcal{B},L^{p}(w)),\quad \frac{1}{\tau} = \alpha + \frac{1}{p}.\end{eqnarray} \end{corollary}Finally we point out that the inclusions in Corollary \ref{coro-CarLpesAppr} can be described in terms of  weighted Besov spaces (\cite{Rou1,Rou2}), namely\begin{eqnarray}\label{defnWeiBesovSpac}\dot{B}_{p,q}^{\alpha}(w) = \{f\in\mathscr{S}'/\mathscr{P}: (2^{kd}\|\varphi_{k}*f\|_{L^{p}(w)})_{k}\in \ell_{q}(\mathbb{Z})\}.\end{eqnarray} (See the definition of $\varphi_{k}$ in section 5). \begin{theorem}\label{thmCaraLppeso} Let $\gamma > 0, \ 1<p<\infty,$ and $w\in A_{p}$ a weight on $\mathbb{R}^{d}.$ Suppose that $\mathcal{B}=\Psi$ is a family of d-dimensional  Lemari\'e-Meyer wavelets or a family of d-dimensional compactly supported Daubechies $D_{N}$ wavelets with $N$ sufficiently large. Then \begin{eqnarray}\label{thmCaraLppeso1}\mathcal{A}_{\tau}^{\frac{\gamma}{d}}(\Psi,L^{p}(w)) = \dot{B}_{\tau,\tau}^{\gamma}(w^{\frac{\tau}{p}})\quad whenever\quad \frac{1}{\tau} = \frac{\gamma}{d} + \frac{1}{p}.\end{eqnarray}  \end{theorem}

The organization of this article is as follows. Basic
facts concerning weights, wavelet bases and greedy bases are given
in section \ref{secPrelim}. Section \ref{SecDemocFunOrl} is devoted to prove Theorem \ref{thmPrinc}. Jackson and Bernstein type  estimates, as well as the inclusions described in Theorem \ref{thmInclAprOrlpesoSpa} are proved in section \ref{secIncluAppr}.  Corollary \ref{coro-CarLpesAppr} and Theorem \ref{thmCaraLppeso}  are proved in section  \ref{secAppSpcLppes}.

\section{Preliminaries}\label{secPrelim}
 \subsection{Basics  in weighted Orlicz spaces}
 In this subsection we recall some basic facts about weights,
weighted  Orlicz spaces and wavelet bases on weighted  Orlicz
spaces, referring to \cite{Benett, Rubio}  for a complete account on these topics.
By a weight on a given measure space, we shall always mean a
non-negative locally integrable function $w$ with values  in $[0,\infty)$  a.e. Let $w(x)$ be
a weight on $\mathbb{R}^{d}$, and for a measurable $Q\subset\mathbb{R}^{d}$ write $w(Q) =
\int_{Q}w(x)dx$. We say that $w \in A_{p}= A_{p}(\mathbb{R}^{d})$, ($1 < p < \infty$) if
there exists a constant $C_{w}$ such that
\begin{eqnarray}\label{condicioAp}
\Big(\frac{1}{|Q|}\int_{Q}w(x)dx\Big)\Big(\frac{1}{|Q|}\int_{Q}w(x)^{-\frac{1}{p-1}}dx\Big)^{p-1}\leq
C_{w},\end{eqnarray} for all  $Q \subset \mathbb{R}^{d}$, where
$|Q|$  denote   the usual  Lebesgue measure of \   $Q$. The
condition $A_{1}$ can be viewed as limiting  case of the condition
$A_{p}$ for $p\downarrow 1$, i.e., (\ref{condicioAp}) is viewed as
\begin{eqnarray}\label{condiciona1}\Big(\frac{1}{|Q|}\int_{Q}w(x)dx\Big)ess_{Q}\sup(w^{-1})\leq C_{w}.\end{eqnarray}
 If $w \in A_{p}$ \ for some $p \in
[1,\infty)$, then there exist $C_{w}^1, C_{w}^2 > 0$ and $\delta > 0$ such
that
\begin{eqnarray}\label{ainfinitocondicion}C_{w}^1\Big(\frac{|A|}{|Q|}\Big)^{p}\leq  \frac{w(A)}{w(Q)}
\leq C_{w}^2 \Big(\frac{|A|}{|Q|}\Big)^{\delta}\end{eqnarray} for all
subsets
 $A\subset Q$.  (For the left hand  inequality take $f=\chi_A$ in part b) of Theorem 2.1,
 Chapter IV, of \cite{Rubio}; for the right hand inequality see Theorem 2.9,  Chapter IV, of \cite{Rubio}).

 A Young function is a convex non-decreasing function $\Phi : [0,\infty)\longrightarrow [0,\infty]$ so that $\lim_{t\longrightarrow +\infty}\Phi(t) = \infty.$ Throughout this paper we shall assume that $\Phi(0) = 0,$    $\Phi$ is strictly increasing and everywhere finite, so that it is a continuos bijection of $[0,\infty).$ Given such $\Phi$ and $w\in A_{\infty} = \cup_{p\geq 1}A_{p},$ the weighted Orlicz space $L^{\Phi}(w)$ is the class of all measurable functions $f: \mathbb{R}^{d}\longrightarrow \mathbb{C}$ so that $\Phi\Big(\frac{|f(x)|}{\lambda}\Big)\in L^{1}(w)$ for some $\lambda > 0.$  The space $L^{\Phi}(w)$ becomes a weighted  rearrangement invariant Banach function space when endowed with the corresponding Luxemburg norm \begin{eqnarray}\label{LuxNormOrlPes}\|f\|_{L^{\Phi}(w)} = \inf\Big\{\lambda > 0:\int_{\mathbb{R}^{d}}\Phi\Big(\frac{|f(x)|}{\lambda}\Big)w(x)dx\leq 1 \Big\}\end{eqnarray} It is not difficult to prove that if $E\subset \mathbb{R}^{d}$ is any measurable set  \begin{eqnarray}\label{fundFuncOrlPeso}\|\chi_{E}\|_{L^{\Phi}(w)} = \frac{1}{\Phi^{-1}(\frac{1}{w(E)})}.\end{eqnarray} The function $\varphi(t) = \frac{1}{\Phi^{-1}(\frac{1}{t})}, \ 0<t<\infty,$ satisfies $\varphi(t) = \|\chi_{E}\|_{L^{\Phi}(w)}$ for any measurable set $E\subset \mathbb{R}^{d}$ such that $w(E) = t,$ and it is called the \textit{fundamental function of $L^{\Phi}(w)$}.

 The Boyd indices of the weighted  Orlicz space $L^{\Phi}(w)$ can be computed directly from the Young function $\Phi$ or from the fundamental function $\varphi.$ Set

\begin{eqnarray}\label{diltFunc}h_{\varphi}^{+}(t) = \sup_{s>0}\frac{\varphi(st)}{\varphi(s)},\quad 0<t<\infty.\end{eqnarray} Then, the lower and upper  Boyd indices $i_{\varphi}$ and $I_{\varphi}$ of $L^{\Phi}(w)$ are given by \begin{eqnarray}\label{LboydInd}i_{\varphi} = \lim_{t\longrightarrow 0}\frac{\log h_{\varphi}^{+}(t)}{\log t} = \sup_{0<t\leq 1}\frac{\log h_{\varphi}^{+}(t)}{\log t} \end{eqnarray}and \begin{eqnarray}\label{UboydIn} I_{\varphi} = \lim_{t\longrightarrow\infty}\frac{\log h_{\varphi}^{+}(t)}{\log t} = \inf_{1<t<\infty}\frac{\log h_{\varphi}^{+}(t)}{\log t}\end{eqnarray}respectively (see \cite{Benett}, p. 277 or \cite{KPS}, p. 54). It is known that $0\leq i_{\varphi}\leq I_{\varphi}\leq 1$ (see Proposition 5.15 of \cite{Benett}, p. 149). Assuming further that $i_{\varphi} > 0$  it follows that\begin{eqnarray}\label{desigFunDilt}\varphi(st)\leq C_{\epsilon}\max\{s^{i_{\varphi}-\epsilon},s^{I_{\varphi}+\epsilon}\}\varphi(t),\quad s,t>0\end{eqnarray}and \begin{eqnarray}\label{desigFunDilt1}\varphi(st)\geq C_{\epsilon}\min\{s^{i_{\varphi}-\epsilon},s^{I_{\varphi}+\epsilon}\}\varphi(t),\quad s,t>0\end{eqnarray} for every $\epsilon > 0$ and some constant $C_{\epsilon} > 0$ (see \cite{Krebec}, p. 3). In this paper we shall only consider weighted Orlicz spaces with non trivial Boyd indices, that is $0<i_{\varphi}\leq I_{\varphi}<1.$

\begin{example}\label{ExamLppesos} When $\Phi(t) = t^{p},\ 1\leq p<\infty,$ then $L^{\Phi}(w) = L^{p}(w)$ and $\varphi(t) = t^{\frac{1}{p}}.$  Hence, $h_{\varphi}^{+}(t) = t^{\frac{1}{p}},$ which implies $i_{\varphi} = I_{\varphi} = \frac{1}{p}.$  \end{example}
\subsection{Wavelet bases and
 weighted   Orlicz spaces.} Let $\mathcal{D}= \{ Q_{j,k} = 2^{-j}([0,1)^{d} + k): j \in \mathbb{Z},k\in\mathbb{Z}^{d}\}
 $ denote the set of all dyadic cubes in $\mathbb{R}^{d}$. We say
 that a finite collection of functions $\{\psi^{1},\ldots,\psi^{L}\}\subset
 L^{2}(\mathbb{R}^{d})$ is an \textit{orthonormal wavelet family}
 if the system \begin{eqnarray}\label{OrthonormalwaveletFamily}
\Big\{\psi_{Q_{j,k}}^{l}(x)= 2^{\frac{jd}{2}}\psi^{l}(2^{j}x-k):j \in \mathbb{Z},k\in\mathbb{R}^{d},l=1,\ldots,L\Big\},\end{eqnarray}
 forms an orthonormal basis of $L^{2}(\mathbb{R}^{d})$. We will say that the wavelet family is $admissible$ if for all $1<p<\infty,$ \begin{eqnarray}\label{caractSquFunLp}\|S_{\psi}f(.)\|_{L^{p}(\mathbb{R}^{d})}\approx\|f(.)\|_{L^{p}(\mathbb{R}^{d})},\end{eqnarray} where \begin{eqnarray}\label{squaFunc1}S_{\psi}f(.) = \Big(\sum_{l=1}^{L}\sum_{I\in\mathcal{D}}|\langle f(.),\psi_{I}(.)\rangle|^{2}\chi_{I}(.)|I|^{-1}\Big)^{\frac{1}{2}}. \end{eqnarray}This implies that wavelet admissible bases are unconditional in $L^{p}(\mathbb{R}^{d}),\ 1<p<\infty.$
 The reader can consult \cite{Dau,Weiss,LemarieMeyer, Me}, for constructions, examples, and properties of orthonormal wavelets.
 Admissible wavelets include the $d$-dimensional Haar system, wavelet arising from multiresolution analysis
 (see \cite{Me}, p. 22),
  wavelets belonging to the regularity class $\mathcal{R}^{0}$
  (as defined in \cite{Weiss}, p. 64 for $d=1$), compactly support wavelets (see \cite{Dau}),   and actually any orthonormal wavelet in
  $L^{2}(\mathbb{R}^{d})$ with mild decay conditions (see \cite{Wo1, Po}).

  In the following result we prove that wavelet admissible basis are also unconditional basis of weighted Orlicz spaces $L^{\Phi}(w),$ for appropriate $w,$  since the norm can be characterize in terms of a square function. Without loss of generality we assume $L=1$ in the rest of this work.

 \begin{theorem}\label{thmcaraOrlPes} Let $L^{\Phi}(w)$ be  a weighted
Orlicz space, with the Boyd indices
satisfying  $0 <i_{\varphi}\leq
I_{\varphi}<1,$  and  $\mathcal{B} = \{\psi_{Q}: Q\in\mathcal{D}\}$ an admissible wavelet basis.  Then, if $w\in A_{p^{\Phi}}(\mathbb{R}^{d}),$ where $p^{\Phi} = \frac{1}{I_{\varphi}},$ we have
\begin{eqnarray}\label{thmcaraOrlpes1}
\|f(.)\|_{L^{\Phi}(w)}\simeq
\|S_{\psi}f(.)\|_{L^{\Phi}(w)}, \quad for \ all \ f   \in L^{\Phi}(w).\end{eqnarray}
\end{theorem}For the
proof we  shall use the following extrapolation theorem adapted to our situation.

\begin{theorem}(\cite{MatellUribePerez})\label{extrap-Theor} Let   $\mathcal{F}$ be a  family of couples of measurable   non-negative functions $(f,g).$ Suppose that for some $1\leq p_{0}<\infty,$ and every weight $w\in A_{p_{0}}(\mathbb{R}^{d}),$  \begin{eqnarray}\label{extrap-Theor1} \int_{\mathbb{R}^{d}}f(x)^{p_{0}}w(x)dx\leq C\int_{\mathbb{R}^{d}}g(x)^{p_{0}}w(x)dx,\quad for\quad all\quad (f,g)\in \mathcal{F}.\end{eqnarray} Then, if $L^{\Phi}(w)$ is a weighted Orlicz  space such that the Boyd indices satisfies, $0< i_{\varphi}\leq I_{\varphi}< 1$  and $w\in A_{p^{\Phi}}(\mathbb{R}^{d}),\ p^{\Phi} = \frac{1}{I_{\varphi}},$ we have \begin{eqnarray}\label{extrap-Theor2}
\|f\|_{L^{\Phi}(w)} \leq C\|g\|_{L^{\Phi}(w)},\quad for \quad all \quad
(f,g)\in\mathcal{F}.\end{eqnarray}
\end{theorem} \begin{proof}( \textbf{of Theorem \ref{thmcaraOrlPes}}) It is proved in  \cite{CUERVAMARTELL} (see also  \cite{Aimar})  that
  \begin{eqnarray}\label{CaracterizacionLppesos1} \|f\|_{L^{p}(w)}
  \simeq \|S_{\psi}(f)\|_{L^{p}(w)}, \end{eqnarray}
   for  all $1< p < \infty$ and  $w \in A_{p}.$
We consider the family $\mathcal{F} = \{(|f|,S_{\psi}(f)): \ S_{\psi}(f)\in
L^{p}(w).$   From the equivalence
(\ref{CaracterizacionLppesos1}), we obtain $$
\int_{\mathbb{R}^{d}}|f(x)|^{p}w(x)dx\leq
C_{1}\int_{\mathbb{R}^{d}}|S_{\psi}(f)|^{p}w(x)dx$$ for all
$1<p<\infty$ and $w\in A_{p}(\mathbb{R}^{d}).$ Then, by Theorem \ref{extrap-Theor}
we obtain $$ \|f\|_{L^{\Phi}(w)}\leq
C_{1}\|S_{\psi}(f)\|_{L^{\Phi}(w)}$$  when
 $w\in A_{p^{\Phi}}.$  The other inequality is proved similarly taking
$\mathcal{F}=\{(S_{\psi}(f),|f|), f\in L^{p}(w)\}.$ \end{proof}
   \subsection{Greedy basis and
 democracy.} We defined in the introduction the notion of {\bf greedy
 basis} in a quasi-normed Banach space $(\mathbb{B},
 \|\cdot\|_{\mathbb{B}})$. We also mentioned the result of Konyagin and
 Temlyakov \cite{Konyagin} characterizing greedy bases as those
 which are unconditional and {\bf democratic}. For simplicity,   given a
 basis $\mathcal{B} = \{e_{j}:j\in \mathbb{N}\}$ in $\mathbb{B}$ we
 shall denote the (normalized)   characteristic   function of a finite set of
 indices \ $\Gamma \in \mathbb{N}$ \  by
 \begin{eqnarray*}\widetilde{1}_{\Gamma} =
 \widetilde{1}_{\Gamma}^{\mathbb{B},\mathcal{B}} =
 \sum_{j\in\Gamma}\frac{e_{j}}{\|e_{j}\|_{\mathbb{B}}}.\end{eqnarray*}
  The basis $\mathcal{B}$ is {\bf democratic} in $\mathbb{B}$ if there exists
 $C \geq 1$ such that
 \begin{eqnarray}\label{conditionDemocraciaFuncionCaract}
 \|\widetilde{1}_{\Gamma}\|_{\mathbb{B}}\leq
 C\|\widetilde{1}_{\Gamma'}\|_{\mathbb{B}}\end{eqnarray} for all
 finite sets of indices $\Gamma, \Gamma' \subset\mathbb{N}$ with $\#\Gamma =
 \#\Gamma'$ (the symbol $\#\Gamma$ denotes the cardinality of the set $\Gamma$). Quite often one can show democracy by finding a function $h: \mathbb{N}\longrightarrow \mathbb{R}^{+}$ for which \begin{eqnarray}\label{condDem}\frac{1}{C}h(\#\Gamma)\leq \|\widetilde{1}_{\Gamma}\|_{\mathbb{B}}\leq C h(\#\Gamma), \quad \  \forall  \ \Gamma \subset \mathbb{N}, \ finite.\end{eqnarray} In the case of wavelet bases, many classical function and distribution spaces satisfy (\ref{condDem}) with $h(\#\Gamma) = (\#\Gamma)^{\frac{1}{p}}.$ Indeed, this is the situation for the Lebesgue spaces $L^{p}(\mathbb{R}^{d})$ when $1<p<\infty$, for the Hardy spaces  $H^{p}(\mathbb{R}^{d}), \ 0< p\leq 1$ and for the  Sobolev spaces $\dot{W}^{s,p}(\mathbb{R}^{d}), \ 1<p<\infty$ (see \cite{HJLY}), and more generally for the family of Triebel-Lizorkin spaces $\dot{F}_{p,r}^{s}(\mathbb{R}^{d})$ with $0<p<\infty,\ s\in \mathbb{R}, \ 0<r\leq\infty$ (under the usual  smoothness assumptions, and with the standard modification of the basis in the case of inhomogeneous spaces; see \cite{GustavoSharpJackson}). Thus, wavelet bases are democratic and hence greedy in all these spaces.

 The Haar system is not greedy in rearrangement invariant spaces defined in $[0,1]$ other than $L^{p}[0,1]$ (see \cite{WoH}). Moreover, wavelet bases  are not democratic in other classical spaces, such as $BMO,$ the Besov spaces $\dot{B}_{p,q}^{\alpha}(\mathbb{R}^{d})$ with $p\neq q,$ Orlicz spaces $L^{\Phi}(\mathbb{R}^{d})$ distinct from $L^{p}(\mathbb{R}^{d}),$ and as we shall see below, weighted Orlicz spaces $L^{\Phi}(w)$ distinct from $L^{p}(w).$\begin{definition}\label{defnRigDemFun} Let $\mathcal{B}$ be a collection of elements in a quasi-Banach space $\mathbb{B}.$ The \textit{right-democracy function associated with $\mathcal{B}$} is defined by\begin{eqnarray}\label{funciDemoSupInf}h_{r}(N;\SB,\B) =
\sup_{Card (\Gamma) = N}
\|\widetilde{1}_{\Gamma}\|_{\SB};
\end{eqnarray} analogously, the \textit{left-democracy function associated with $\mathcal{B}$} is defined by \begin{eqnarray}\label{left-Democ-Func} h_{l}(N;\SB,\B)= \inf_{Card(\Gamma) =
N}\|\widetilde{1}_{\Gamma}\|_{\SB}\end{eqnarray}\end{definition}

Observe that a basis $\mathcal{B}$ is democratic  in $\mathbb{B}$ if and only if, $h_{r}(N;\mathbb{B},\mathcal{B})\leq C h_{l}(N;\mathbb{B},\mathcal{B})$ for all $N\geq 1$ and some $C > 0.$

We want to  show  that, in general, admissible wavelet bases are not democratic in weighted Orlicz spaces. In order to do so one needs to estimate $\|\widetilde{1}_{\Gamma}\|_{L^{\Phi}(w)}$ in terms of $\#\Gamma.$ This can be done when $\Gamma$ is a collection of  pairwise disjoint dyadic cubes $\{Q_{j}\}_{j=1}^{N},$ such that $w(Q_{j})\approx \tau,$  for any $\tau > 0.$

We state and prove the following   results.\begin{lemma}\label{auxLemm1}Let $w\in A_{\infty}(\mathbb{R}^{d})$ be a weight. If  $\{Q_{k}\}_{k=-\infty}^{\infty}$ is a family of dyadic cubes such that $Q_{k}\subset Q_{k+1}$ and $|Q_{k+1}| = 2^{d}|Q_{k}|$ for all $k\in\mathbb{Z}^{+},$ then \begin{eqnarray}\label{auxLemm12}\lim_{k\longrightarrow\infty}w(Q_{k})=\infty \quad \ \and \quad \lim_{k\longrightarrow -\infty}w(Q_{k}) = 0.\end{eqnarray}\end{lemma}\begin{proof}\label{proof-auxLem1} Because $w\in A_{\infty},$  if $k\geq 0,$ by (\ref{ainfinitocondicion}) we obtain $$\frac{w(Q_{0})}{w(Q_{k})}\leq C_{w}^{2}\Big(\frac{|Q_{0}|}{|Q_{k}|}\Big)^{\delta} = C_{w}^{2}\Big(\frac{1}{2^{kd}}\Big)^{\delta}.$$ Then, $w(Q_{k})\geq (C_{w}^{2})^{-1}2^{kd\delta}w(Q_{0})$ and $\lim_{k\longrightarrow\infty}w(Q_{k})=\infty.$ On the other hand, if $k\leq 0,$ by (\ref{ainfinitocondicion}) we obtain $$\frac{w(Q_{k})}{w(Q_{0})}\leq C_{w}^{2}\Big(\frac{|Q_{k}|}{|Q_{0}|}\Big)^{\delta} = C_{w}^{2}2^{kd\delta}.$$ Then, $w(Q_{k})\leq C_{w}^{2}2^{kd\delta}w(Q_{0})$ and $\lim_{k\longrightarrow -\infty}w(Q_{k}) = 0.$\end{proof}\begin{lemma}\label{auxLemm2}Let $w\in A_{\infty}(\mathbb{R}^{d})$ be a weight. Given $\tau > 0$ there exists a pairwise disjoint sequence of cubes $\{R_{j}\}_{j=1}^{\infty}\subset\mathcal{D}$ such that $$C\tau< w(R_{j})\leq\tau$$ where $C> 0$ is a constant depending only on $w.$ \end{lemma}\begin{proof}\label{proof-auxLem2} Let $Q_{k} = [0,2^{k})^{d}, \ k\in \mathbb{Z}.$ By lemma \ref{auxLemm1} there exists $k_{1}\in\mathbb{Z}$ such that \begin{eqnarray}\label{proof-auxLem21}w(Q_{k_{1}})\leq \tau < w(Q_{k_{1}+1}).\end{eqnarray} Choose $R_{1} = Q_{k_{1}}.$ We have $w(R_{1}) = w(Q_{k_{1}})\leq \tau.$ On the other hand, by (\ref{ainfinitocondicion}), we obtain $$\frac{w(Q_{k_{1}})}{w(Q_{k_{1}+1})}\geq C_{w}^{1}\Big(\frac{|Q_{k_{1}}|}{|Q_{k_{1}+1}|}\Big)^{p} = C_{w}^{1}2^{-dp},$$ so that $$w(R_{1})=w(Q_{k_{1}})\geq C_{w}^{1}2^{-dp}w(Q_{k_{1}+1})> C_{w}^{1}2^{-dp}\tau.$$ Thus, we can take $C = C_{w}^{1}2^{-dp}.$

Suppose we have chosen disjoint cubes $R_{1}, R_{2},\ldots, R_{m-1}$ such that $C\tau< w(R_{j})\leq \tau$ for all $j=1,2,\ldots,m-1.$ Without loss of generality we can assume that all the $R_{j}$ are contained in the positive cone of  $\mathbb{R}^{d},$ that is, the set of points of $\mathbb{R}^{d}$ with non-negative coordinates.

Choose $Q_{0} = 2^{k_{m}}[0,1)^{d},\ k_{m}\in\mathbb{Z},$ such that $R_{j}\subset Q_{0}$ for all $j=1,2,\ldots,m-1.$ consider the increasing family of dyadic cubes given by $Q_{k} = 2^{k_{m}+k}[0,1)^{d},\ k=0,1,2,\ldots.$ Let $\widetilde{Q}_{k},\ k=1,2,\ldots,$ be a dyadic cube contained in $Q_{k}$ such that $|\widetilde{Q}_{k}| = \frac{|Q_{k}|}{2^{d}}$ and $\widetilde{Q}_{k}\cap Q_{k-1} = \emptyset.$  If $w(\widetilde{Q}_{k})\leq \tau$ for all $k= 1,2,3,\ldots$ by (\ref{ainfinitocondicion}) we obtain $$\frac{w(\widetilde{Q}_{k})}{w(Q_{k})}\geq C_{w}^{1}\Big(\frac{|\widetilde{Q}_{k}|}{|Q_{k}|}\Big)^{p} = C_{w}^{1}2^{-dp}.$$ Thus, $w(Q_{k})\leq (C_{w}^{1})^{-1}2^{dp}\tau$ for all $k=1,2,\ldots$ contradicting lemma \ref{auxLemm1}. Thus, there exists $k_{m}^{0}\in\mathbb{Z}$ such that $w(\widetilde{Q}_{k_{m}^{0}})> \tau.$ Consider a family of descendants of the dyadic cube $\widetilde{Q}_{k_{m}^{0}}.$ By lemma \ref{auxLemm1}, there exists $\widetilde{Q}_{k_{m}}, \ \widetilde{\widetilde{Q}}_{k_{m}}\in\mathcal{D}$ such that \begin{eqnarray}\label{proof-auxLem22}w(\widetilde{\widetilde{Q}}_{k_{m}})\leq \tau < w(\widetilde{Q}_{k_{m}})\end{eqnarray}and $|\widetilde{\widetilde{Q}}_{k_{m}}| = \frac{|\widetilde{Q}_{k_{m}}|}{2^{d}}.$  Choose $R_{m} = \widetilde{\widetilde{Q}}_{k_{m}}.$ Since (\ref{proof-auxLem22}) is the same relation as (\ref{proof-auxLem21}) it follows that $$ C_{w}^{1}2^{-dp}\tau< w(R_{m})\leq \tau.$$ Observe that $R_{m}$ has been chosen in the positive  cone of $\mathbb{R}^{d}$ and is disjoint to $R_{1},\ldots,R_{m-1}.$
  \end{proof}
\begin{proposition}\label{PropDisCub} Let
$L^{\Phi}(w)$ be a weighted Orlicz space with Boyd
indices $0< i_{\varphi} \leq I_{\varphi} < 1,\ w\in A_{p^{\Phi}}$ a weight in $\mathbb{R}^{d}$, and
let $ \mathcal{B}= \{ \psi_{Q} : Q \in \mathcal{D}\}$ be an
admissible wavelet basis. \\i) If  $\Gamma = \{Q_{1},\ldots,Q_{N}\}
\subset\mathcal{D}$ is a pairwise disjoint family then
\begin{eqnarray}\label{PropDisCub1}
\|\widetilde{1}_{\Gamma}\|_{L^{\Phi}(w)}\approx \Big\|\sum_{Q\in\Gamma}\frac{\chi_{Q}(.)}{\varphi(w(Q))}\Big\|_{L^{\Phi}(w)}.\end{eqnarray} ii)  Moreover, for any $\tau > 0,$ there exist a family of disjoint cubes \   $\Gamma = \{R_{1}, R_{2},\ldots,R_{N} \} \\ \subset\mathcal{D},$ such that   \begin{eqnarray}\label{prpDiCub2} \|\widetilde{1}_{\Gamma}\|_{L^{\Phi}}\approx \frac{\varphi(N\tau)}{\varphi(\tau)}\end{eqnarray}
\end{proposition}

\begin{proof} $i)$ For a single element of the basis $\mathcal{B}$ we
have, by (\ref{thmcaraOrlpes1}) that
\begin{eqnarray}\label{OrlicznormSingleelement}
\|\psi_{Q}\|_{L^{\Phi}(w)}\simeq\Big\|\Big(\frac{\chi_{Q}(.)}{|Q|}\Big)^{\frac{1}{2}}\Big\|_{L^{\Phi}(w)}
=
\frac{\|\chi_{Q}(.)\|_{L^{\Phi}(w)}}{|Q|^{\frac{1}{2}}}=\frac{\varphi(w(Q))}{|Q|^{\frac{1}{2}}}.\end{eqnarray}
By  (\ref{thmcaraOrlpes1}) again
\begin{eqnarray*}\|\widetilde{1}_{\Gamma}\|_{L^{\Phi}(w)}
&\simeq&
 \Big\|\Big(\sum_{Q\in\Gamma}\frac{1}{\|\psi_{Q}\|_{L^{\Phi}(w)}^{2}}\chi_{Q}|Q|^{-1}\Big)^{\frac{1}{2}}\Big\|_{\mathrm{L}^{\Phi}(w)}\approx\Big\|\Big(\sum_{Q\in\Gamma}\frac{\chi_{Q}}{\varphi(w(Q))^{2}}\Big)^{\frac{1}{2}}\Big\|_{\mathrm{L}^{\Phi}(w)}
={} \nonumber\\ {} && =
\Big\|\sum_{Q\in\Gamma}\frac{\chi_{Q}}{\varphi(w(Q))}\Big\|_{L^{\Phi}(w)},
\end{eqnarray*} where in the last equality we have used that the cubes in $\Gamma$ are pairwise disjoint.

 $ii)$ The existence of the family  $\Gamma = \{ R_{1},R_{2},\ldots,R_{N}\}\subset\mathcal{D}$ is proved in Lemma \ref{auxLemm2} where it is shown that $w(R_{j})\approx \tau, \ j=1,2,\ldots,N.$ In this situation \begin{eqnarray}\label{proof-PropDicCub}\|\widetilde{1}_{\Gamma}\|_{L^{\Phi}(w)}&=& \Big\|\sum_{j=1}^{N}\frac{\chi_{R_{j}}(.)}{\varphi(w(R_{j}))}\Big\|_{L^{\Phi}(w)}\approx \frac{1}{\varphi(\tau)}\Big\|\bigcup_{j=1}^{N}\chi_{R_{j}}\Big\|_{L^{\Phi}(w)}{}\nonumber\\{}&&=\frac{1}{\varphi(\tau)}\varphi(w\Big(\bigcup_{j=1}^{N}R_{j}\Big))\approx \frac{\varphi(N\tau)}{\varphi(\tau)}.\end{eqnarray}
 \end{proof}\begin{Remark}\label{remarkPropDisCub} It follows from part $ii)$ of Proposition \ref{PropDisCub} that for admissible wavelet basis $\mathcal{B}$ $$h_{r}(N;L^{\Phi}(w),\mathcal{B})\gtrsim \sup_{\tau > 0}\frac{\varphi(N\tau)}{\varphi(\tau)} = h_{\varphi}^{+}(N)$$ and $$h_{l}(N;L^{\Phi}(w),\mathcal{B})\lesssim \inf_{\tau>0}\frac{\varphi(N\tau)}{\varphi(\tau)} = h_{\varphi}^{-}(N).$$\end{Remark} Thus, if $h_{\varphi}^{+}(N)$ and $h_{\varphi}^{-}(N)$ are not comparable for $N\geq 1$ it follows that admissible wavelet bases are non democratic in weighted Orlicz spaces. On the other hand, if wavelet admissible bases are democratic in $L^{\Phi}(w),$ $h_{\varphi}^{+}(N)\leq h_{\varphi}^{-}(N),$ and Lemma 5.2 in \cite{GHM} shows that $\varphi(t)\approx t^{\alpha}$ for some $\alpha\in(0,1);$ thus, the only democratic weighted Orlicz spaces are the spaces $L^{p}(w)$ for some $p = \frac{1}{\alpha}\in (1,\infty).$

 \section{Left and Right democracy functions for weighted  Orlicz spaces}\label{SecDemocFunOrl}
Our main theorem  in this section shows that $h_{r}(N; L^{\Phi}(w),\mathcal{B})\lesssim h_{\varphi}^{+}(N)$ and\\ $h_{l}(N; L^{\Phi},\mathcal{B}) \gtrsim h_{\varphi}^{-}(N)$ (see theorem \ref{thmdemFuncOrlpes} below) giving us   together with remark \ref{remarkPropDisCub}  a complete description (up to multiplicative constants) of the left and right democracy functions of wavelet basis on weighted  Orlicz spaces.

\begin{theorem}\label{thmdemFuncOrlpes} Let  \ $L^{\Phi}(w),$ be a  weighted Orlicz
space with   Boyd indices satisfying $0<i_{\varphi}\leq I_{\varphi}<1$ , $w\in A_{p^{\Phi}}$ a weight in $\mathbb{R}^{d},$  and let $\mathcal{B}=
\{\psi_{Q}: Q\in\mathcal{D}\}$ be an admissible wavelet basis. Then for all $\Gamma\subset\mathcal{D}$
\begin{eqnarray}\label{democracyinOrlicz1}
h_{\varphi}^{-}(\#\Gamma)\lesssim \|\widetilde{1}_{\Gamma}\|_{L^{\Phi}(w)}\lesssim h_{\varphi}^{+}(\#\Gamma) .\end{eqnarray} This, together with Remark \ref{remarkPropDisCub} gives $$h_{r}(N;\mathcal{B},L^{\Phi}(w))\approx h_{\varphi}^{+}(N)\quad and \quad h_{l}(N;\mathcal{B},L^{\Phi}(w))\approx h_{\varphi}^{-}(N),$$ which is Theorem \ref{thmPrinc}.\end{theorem}

The rest of this section is devoted to prove Theorem  \ref{thmdemFuncOrlpes}. We first present a very simple argument for the case of pairwise disjoint  cubes.
 \subsection{Proof of Theorem \ref{thmdemFuncOrlpes}: The case of disjoint cubes }
We assume first that $\Gamma = \{Q_{1},\ldots,Q_{N}\}$ consists of pairwise dijoint cubes.
  Let
$\lambda = h_{\varphi}^{+}(N)$, so that $\varphi(N\omega(Q)) \leq
\lambda\varphi(\omega(Q)), \ \forall\  Q \in \Gamma$. Therefore,
since the elements of $\Gamma$ are disjoint, and $\Phi$ is increasing
\setlength\arraycolsep{2pt} \begin{eqnarray*}
&&\int_{\mathbb{R}^{d}}\Phi\Big(\frac{\sum_{j=1}^{N}\frac{\chi_{Q_{j}}(x)}{\varphi(\omega(Q_{j}))}}{\lambda}\Big)w(x)dx
=
\sum_{j=1}^{N}\Phi\Big(\frac{1}{\lambda\varphi(\omega(Q_{j}))}\Big)\omega(Q_{j}){} \nonumber\\ {} &&\leq
\sum_{j=1}^{N}\Phi\Big(\frac{1}{\varphi(N\omega(Q_{j}))}\Big)\omega(Q_{j})
=\sum_{j=1}^{N}\Phi\Big(\Phi^{-1}\Big(\frac{1}{N\omega(Q_{j})}\Big)\Big)\omega(Q_{j})=
1\end{eqnarray*}Then by (\ref{PropDisCub1}) and (\ref{LuxNormOrlPes})
we have \begin{eqnarray*}
\|\widetilde{1}_{\Gamma}\|_{L^{\Phi}(w)}\simeq \Big\|\sum_{j=1}^{N}\frac{\chi_{Q_{j}}(.)}{\varphi(w(Q_{j}))}\Big\|_{L^{\Phi}(w)}\leq
h_{\varphi}^{+}(N).\end{eqnarray*} The lower estimate is obtained in
a similar way.

\subsection{Proof of Theorem \ref{thmdemFuncOrlpes}: The general case}

 In the case of
disjoint cubes just considered we have two important features.
First, Proposition \ref{PropDisCub} allows us to
``linearize'' the square function in
(\ref{thmcaraOrlpes1}). Second, for the
estimates obtained in the previous argument it is crucial that the
sets involved are disjoint. For general families of cubes we are
going to follow the same scheme. First we  ``linearize'' the square
function and we dominate this by an expression involving only
disjoint subsets from the elements of
$\Gamma$.\\

\textbf{Linearization of the square function}.  Given a finite set
$\Gamma \subset\mathcal{D}$, we denote
\begin{eqnarray}\label{squarefuntionPeso} S_{\Gamma}(x) =
\Big(\sum_{Q\in\Gamma}\frac{\chi_{Q}(x)}{\varphi(w(Q))^{2}}\Big)^{\frac{1}{2}},\end{eqnarray}
so that by (\ref{thmcaraOrlpes1}) and
(\ref{OrlicznormSingleelement}), we have
$\|\widetilde{1}_{\Gamma}\|_{L^{\Phi}(w)} \simeq
\|S_{\Gamma}(\cdot)\|_{L^{\Phi}(w)}$. For every $x \in
\bigcup_{Q\in\Gamma}Q$, we define $Q_{x}$ as the smallest (hence
unique) cube in $\Gamma$ containing $x$. It is clear that
\begin{eqnarray}\label{linearization1} S_{\Gamma}(x) \geq
\frac{\chi_{Q_{x}}(x)}{\varphi(w(Q_{x}))}, \quad \forall \ x \in
\bigcup_{Q\in \Gamma}Q,\end{eqnarray} since the left hand side
contains at least the cube  $Q_{x}$ (and possible more). We now
show that the reverse inequality holds. Indeed, if we enlarge the
sum to include all dyadic cubes containing $Q_{x}$ we have
\begin{eqnarray*} S_{\Gamma}(x)^{2} =
\sum_{Q\in\Gamma}\frac{\chi_{Q}(x)}{\varphi(w(Q))^{2}}\leq
\sum_{\underset{Q\in\mathcal{D}}{Q\supset
Q_{x}}}\frac{1}{\varphi(w(Q))^{2}}\leq
\sum_{j=0}^{\infty}\frac{1}{\varphi(w(Q_{x}^{j}))^{2}},\end{eqnarray*} where $Q_{x}^{j}$ denotes the unique cube of measure $2^{jd}|Q_{x}|$ containing $Q_{x}.$  Now since  $Q_{x} = Q_{x}^{0}
\subset Q_{x}^{1} \subset Q_{x}^{2}\subset\ldots$ we can use  (\ref{ainfinitocondicion}) to obtain
\begin{eqnarray*}\frac{w(Q_{x})}{w(Q_{x}^{j})}\leq
C_{w}^{2}\Big(\frac{|Q_{x}|}{|Q_{x}^{j}|}\Big)^{\delta} =
C_{w}^{2}2^{-jd\delta}.\end{eqnarray*} Hence,
\begin{eqnarray*}w(Q_{x}^{j})\geq
(C_{w}^{2})^{-1}w(Q_{x})2^{jd\delta},\end{eqnarray*} and $$\varphi(w(Q_{x}^{j}))\geq \varphi((C_{w}^{2})^{-1}w(Q_{x})2^{jd\delta}).$$ Since $i_{\varphi} > 0$, by
(\ref{desigFunDilt1}) we can choose
$0<\epsilon<i_{\varphi}$ and find a $C_{\epsilon}> 0$ such that
$\varphi((C_{w}^{2})^{-1}2^{jd\delta}w(Q_{x})) \geq
C_{\epsilon}((C_{w}^{2})^{-1}2^{jd\delta})^{(i_{\varphi}-\epsilon)})\varphi(w(Q_{x})$.\\

Thus,\begin{eqnarray*}S_{\Gamma}(x)^{2}\leq
\frac{(C_{w}^{2})^{(i_{\varphi}-\epsilon)}C_{\epsilon}}{(\varphi(w(Q_{x})))^{2}}\sum_{j=0}^{\infty}2^{-jd\delta(i_{\varphi}
-
\epsilon)}\lesssim\frac{\chi_{Q_{x}}(x)}{(\varphi(w(Q_{x})))^{2}}.\end{eqnarray*}
This and (\ref{linearization1}) show that
\begin{eqnarray}\label{linearization2}S_{\Gamma}(x) \simeq
\frac{\chi_{Q_{x}}(x)}{\varphi(w(Q_{x}))}.\end{eqnarray} Observe
from (\ref{linearization2}) that $S_{\Gamma}(x)\simeq
S_{\Gamma_{\min}}(x)$, where $\Gamma_{\min}(x)$ denotes the family
of minimal cubes in $\Gamma$, that is, \begin{eqnarray*}
\Gamma_{\min} =\Big \{Q_{x}: x \in\bigcup_{Q\in\Gamma}Q\Big\}.\end{eqnarray*}
\subsection{Shaded and Lighted Cubes}
Shaded and lighted cubes were introduced in \cite{GHM}. We recall the definitions. Given a fixed
$\Gamma\subset\mathcal{D}$, for any $Q\in\Gamma$ we define the
\textbf{Shade} of $Q$ as the union of all cubes from $\Gamma$
strictly contained in $Q$
\begin{eqnarray*}Shade(Q) = \bigcup\Big\{R : R\in \Gamma,
R\subsetneq Q \Big\}.\end{eqnarray*} We define the \textbf{Light} of $Q$
as  $Light(Q) = Q \backslash
Shade(Q).$  It is clear that $Q \in
\Gamma_{\min},$ if and only if, $Light(Q)\neq \emptyset$, and moreover
\begin{eqnarray*}\bigcup_{Q\in\Gamma}Q =
\bigcup_{Q\in\Gamma_{\min}}Light(Q).\end{eqnarray*} Therefore, by
(\ref{linearization2}) we can write \begin{eqnarray}
\label{linearization3} S_{\Gamma}(x) \simeq \sum_{Q \in
\Gamma_{\min}}\frac{\chi_{Light(Q)}(x)}{\varphi(w(Q))},\end{eqnarray}
where in the last sum there is at most one non-zero term for each $x$. We shall
classify the cubes as shaded if the shade is a big portion of the
cube or lighted if this does not happen. Precisely, a cube $Q \in
\Gamma$ is called \textbf{shaded} if $|Shade(Q)| >
\frac{2^{d}-1}{2^{d}}|Q|$, and we write $\Gamma_{s}$ for the
collection of cubes from $\Gamma$ that are shaded. A cube $Q$ from
$\Gamma$ is called \textbf{lighted} if it is not shaded, that is, if
$|Light(Q)|\geq\frac{1}{2^{d}}|Q|$. We write $\Gamma_{L}$ for the
collection of all cubes from $\Gamma$ that are lighted.\begin{Remark}\label{RcardSets} Observe that
$\Gamma_{L} \subset
\Gamma_{\min}$ and by Lemma 4.3 in \cite{GHM} we have $$\frac{2^{d}-1}{2^{d}}(\#\Gamma)\leq (\#\Gamma_{L})\leq(\#\Gamma_{\min})\leq(\#\Gamma), \quad \forall\ \Gamma\subset\mathcal{D}$$ \end{Remark}  Now we shall
conclude the proof of theorem \ref{thmdemFuncOrlpes}.\begin{proof}( \textbf{of Theorem \ref{thmdemFuncOrlpes})} By
(\ref{thmcaraOrlpes1}) and (\ref{linearization3}) we know that
\begin{eqnarray}\label{carOrlLight}
\|\widetilde{1}_{\Gamma}\|_{L^{\Phi}(w)}\simeq\Big\|\sum_{Q\in\Gamma_{\min}}\frac{\chi_{Light(Q)}(x)}{\varphi(w(Q))}\Big\|_{L^{\Phi}(w)}.\end{eqnarray}
Thus it is enough to estimate the quantity in the right side of
(\ref{carOrlLight}).  We let
$\lambda = h_{\varphi}^{+}(\#\Gamma_{\min})$ so that
$\varphi(w(Q)\#\Gamma_{\min})) \leq \lambda\varphi(w(Q))$ for all $Q \in
\Gamma_{\min}$. Since $\{Light(Q): Q \in\Gamma_{\min}\}$ is a
pairwise disjoint  collection and $\Phi$ is increasing, we
have\begin{eqnarray*}
&&\int_{\mathbb{R}^{d}}\Phi\Big(\frac{\sum_{Q\in\Gamma_{\min}}\frac{\chi_{Light(Q)(x)}}{\varphi(w(Q))}}{\lambda}\Big)w(x)dx
=
\sum_{Q\in\Gamma_{\min}}\Phi\Big(\frac{1}{\lambda\varphi(w(Q))}\Big)w(Light(Q)){}
\nonumber\\ {} && \leq
\sum_{Q\in\Gamma_{\min}}\Phi\Big(\frac{1}{\varphi(w(Q)\#\Gamma_{\min})}\Big)w(Q)
=\sum_{Q \in
\Gamma_{\min}}\Phi\Big(\Phi^{-1}\Big(\frac{1}{w(Q)\#\Gamma_{\min}}\Big)\Big)w(Q)
= 1.\end{eqnarray*} Hence by,
(\ref{carOrlLight}), Remark \ref{RcardSets}
 and since $h_{\varphi}^{+}$ is non
decreasing, we have \begin{eqnarray*}
\|\widetilde{1}_{\Gamma}\|_{L^{\Phi}(w)}\lesssim
h_{\varphi}^{+}(\#\Gamma_{\min})\lesssim
h_{\varphi}^{+}(\#\Gamma).\end{eqnarray*} For the left inequality, by
(\ref{carOrlLight}) and using that
$\Gamma_{L} \subset \Gamma_{\min}$, we can write
\begin{eqnarray*}\setlength\arraycolsep{2pt}
\|\widetilde{1}_{\Gamma}\|_{L^{\Phi}(w)}&\gtrsim&
\Big\|\sum_{Q\in\Gamma_{L}}\frac{\chi_{Light(Q)}(x)}{\varphi(w(Q))}\Big\|_{L^{\Phi}(w)}.
 \end{eqnarray*} Now let $\lambda <
 h_{\varphi}^{-}(2^{-dp}C_{w}^{1}(\#\Gamma_{L}))$ so that $\lambda\varphi(w(Q))<
 \varphi(w(Q)2^{-dp}C_{w}^{1}(\#\Gamma_{L}))$ for all $Q\in\Gamma_{L}$.
 Using (\ref{ainfinitocondicion}), and  since $|Light(Q)|>2^{-d}|Q|$
 for $Q \in \Gamma_{L}$, we deduce, with $p=p^{\Phi},$ that \begin{eqnarray*}
&&\int_{\mathbb{R}^{d}}\Phi\Big(\frac{\sum_{Q\in\Gamma_{L}}\frac{\chi_{Light(Q)}(x)}{\varphi(w(Q))}}{\lambda}\Big)w(x)dx
= \sum_{Q\in
\Gamma_{L}}\Phi\Big(\frac{1}{\lambda\varphi(w(Q))}\Big)w(Light(Q))
{} \nonumber\\ {} > && \sum_{Q \in \Gamma_{L}}\Phi\Big(\frac{1}{\varphi(2^{-dp}C_{w}^{1}w(Q)(\#\Gamma_{L}))}\Big)C_{w}^{1}2^{-dp}w(Q) {}\nonumber\\ {} && = \sum_{Q\in\Gamma_{L}}\Phi\Big(\Phi^{-1}\Big(\frac{1}{2^{-dp}C_{w}^{1}w(Q)(\#\Gamma_{L})}\Big)\Big) 2^{-dp}w(Q)C_{w}^{1} =
1.\end{eqnarray*} Then by (\ref{LuxNormOrlPes}), and  Remark \ref{RcardSets}   we
obtain
\begin{eqnarray*} \|\widetilde{1}_{\Gamma}\|_{L^{\Phi}(w)}\geq
h^{-}_{\varphi}(2^{-dp}C_{w}^{1}(\#\Gamma_{L}))\geq
h_{\varphi}^{-}(C_{w}^{1}(2^{d}-1)2^{-d(p+1)}(\#\Gamma)).\end{eqnarray*}  Now using   (\ref{desigFunDilt1}) it can be shown that \begin{eqnarray*}h_{\varphi}^{-}(C_{w}^{1}(2^{d}-1)2^{-d(p+1)}(\#\Gamma)) \geq
Ch_{\varphi}^{-}(\#\Gamma).
\end{eqnarray*}\end{proof}
If $\Phi(t) = t^{p},$ from Theorem \ref{thmdemFuncOrlpes} and Example \ref{ExamLppesos} we deduce that admissible wavelet bases are democratic in weighted Lebesgue spaces $L^{p}(w).$

\begin{corollary}\label{CorodemocFuncLppeso} Let $\Phi(t) = t^{p},\ 1<p<\infty,$  $w\in A_{p}$ a weight in $\mathbb{R}^{d},$ and $\mathcal{B} = \{\psi_{Q}: Q\in\mathcal{D}\}$ an admissible wavelet basis. Then \begin{eqnarray}\label{democFuncLppeso} h_{r}(N;\mathcal{B},L^{p}(w))\approx h_{l}(N; \mathcal{B},L^{p}(w))\approx N^{\frac{1}{p}}\end{eqnarray}\end{corollary}
\section{Inclusions  for $N$-Term Approximation spaces of $L^{\Phi}(w).$}\label{secIncluAppr}
In this section we investigate Jackson and Bernstein type inequalities and the corresponding inclusions for the $N$-term approximation spaces $\mathcal{A}_{q}^{\alpha}(\mathcal{B},L^{\Phi}(w)),\ \alpha > 0,\ 0<q\leq\infty,\ w\in A_{p^{\Phi}}(\mathbb{R}^{d}),$ where the error of approximation is measured in $L^{\Phi}(w)$(see \ref{ApproximationSpaceNorm}). These inclusions are given in terms of the discrete Lorentz spaces $\Lambda_{\eta}^{q}$ (see  definition and properties of this spaces in subsection  \ref{subsSeqSpc}, bellow).

\subsection{Sequence spaces in $\mathcal{D}$}\label{subsSeqSpc}
We recall the definition of some classical   sequence spaces over the index set $\mathcal{D}$ of all dyadic cubes on $\mathbb{R}^{d}.$ All of them are subspaces of $c_{o}$ and therefore for each sequence $\{s_{Q}\}_{Q\in\mathcal{D}}$ we can find an enumeration of the index set $\mathcal{D}=\{Q_{k}\}_{k=1}^{\infty}$ so that $|s_{Q_{1}}|\geq|s_{Q_{2}}|\geq\ldots$ and in addition $\lim_{k\longrightarrow\infty}s_{Q_{k}} = 0.$ We shall always assume that $\{s_{Q_{k}}\}_{k\geq 1}$ corresponds to such ordering, which coincides with the \textit{non-increasing rearrangement} $\textbf{s}^{\ast}$ of the sequence $\textbf{s}.$

Let $\eta = \{\eta(k)\}_{k\geq 1}$ be a fixed positive increasing sequence so that $\lim_{k\longrightarrow\infty}\eta(k) = \infty$ and $\eta$ is \textit{doubling} (i.e. $\eta(2k)\leq C\eta(k), k\geq 1$). Then, for each $0<q\leq\infty$ we define a \textit{weighted discrete Lorentz space} by $$\Lambda_{\eta}^{q}= \Big\{\textbf{s}\in c_{o}: \|\textbf{s}\|_{\Lambda_{\eta}^{q}} = \Big[\sum_{k\geq 1}(\eta(k)|s_{Q_{k}}|)^{q}\frac{1}{k}\Big]^{\frac{1}{q}}<\infty\Big\}.$$
Note that for $q=\infty$ one writes $\|\textbf{s}\|_{\Lambda_{\eta}^{\infty}} = \sup_{k}\eta(k)|s_{Q_{k}}|.$ These are quasi-Banach rearrangement invariant spaces, which are Banach when $q\geq 1$ and  $\{\frac{\eta(k)^{q}}{k}\}_{k}$ is non-increasing (\cite{CRSoria}, p. 28). When $q = 1$ or $q = \infty$ we shall write, respectively, $\Lambda_{\eta}$ and $\mathbb{M}_{\eta}$ (the latter called Marcinkiewicz space). The particular case $\eta(k) = k^{\frac{1}{\tau}}$ gives  the classical (discrete) Lorentz space $\Lambda_{\eta}^{q} = \ell^{\tau,q}(\mathcal{D}).$ The spaces $\Lambda_{\eta}^{q}$ for general $\eta,$ and in particular, their interpolation properties, have been studied, e.g., in  \cite{CRSoria,Meruci,Person}. In our applications we use the sequences $\{k^{\alpha}h_{\varphi}^{\pm}(k)\}_{k\geq 1},$ for  $\alpha >0,$ which always satisfy the required assumptions.

Given a fixed sequence space $\mathfrak{s}$ as above, we define a new sequence space $\mathfrak{s}(L^{\Phi}(w))$ isomorphic to $\mathfrak{s},$ by $$\mathfrak{s}(L^{\Phi}(w)) = \{f=\sum_{Q\in\mathcal{D}}\langle f,\psi_{Q}\rangle\psi_{Q}\in L^{\Phi}(w): \{\|\langle f,\psi_{Q}\rangle\psi_{Q}\|_{L^{\Phi}(w)}\}_{Q}\in\mathfrak{s}\},$$ with $\|f\|_{\mathfrak{s}(L^{\Phi}(w))} = \Big\|\|\langle f,\psi_{Q_{k}}\rangle\psi_{Q_{k}}\|_{L^{\Phi}(w)}\Big\|_{\mathfrak{s}}.$ Such definitions appear naturally in relation with approximation when the basis is not normalized (see, e.g., \cite{GustavoSharpJackson}).

\subsection{Jackson type inequalities} In order to obtain the left embedding of the inclusions of approximation spaces given in Theorem \ref{thmInclAprOrlpesoSpa}, we start by proving some inequalities of Jackson Type.\begin{proposition}\label{PropJacInqSeqOrlSpa} Let $\Phi$ be a Young function so that $0<i_{\varphi}\leq I_{\varphi}<1,$ $w\in A_{p^{\Phi}}$ a weight in $\mathbb{R}^{d},$ and let $\alpha> 0.$  Let $\mathcal{B}$ be an admissible wavelet basis. Then, there exists $C>0$ such that for every $f\in\mathbb{M}_{k^{\alpha}h_{\varphi}^{+}(k)}(\mathcal{B},L^{\Phi}(w))$ we have \begin{eqnarray}\label{PropJacInqSeqOrlSpa1}\|f-G_{N-1}(f)\|_{L^{\Phi}(w)}\leq CN^{-\alpha}\|f\|_{\mathbb{M}_{k^{\alpha}h_{\varphi}^{+}(k)}(\mathcal{B},L^{\Phi}(w))}, \quad \forall N\geq 1.\end{eqnarray}\end{proposition}\begin{proof}\label{proofPropJacIneq}By the triangle inequality and (\ref{nonincreaRearrang}) we have \begin{eqnarray}\label{proofPropJacIneq1}\|f - G_{N-1}(f)\|_{L^{\Phi}(w)} &=& \Big\|\sum_{k\geq N} \langle f,\psi_{Q_{k}}\rangle \psi_{Q_{k}}\Big\|_{L^{\Phi}(w)}\leq\sum_{j=0}^{\infty}\Big\|\sum_{2^{j}N\leq k< 2^{j+1}N}\langle f,\psi_{Q_{k}}\rangle \psi_{Q_{k}}\Big\|_{L^{\Phi}}{}\nonumber\\{}&&\leq \sum_{j=0}^{\infty}\|\langle f,Q_{2^{j}N}\rangle\psi_{Q_{2^{j}N}}\|_{L^{\Phi}(w)}\Big\|\sum_{2^{j}N\leq k<2^{j+1}N} \frac{\psi_{Q_{k}}}{\|\psi_{Q_{k}}\|_{L^{\Phi}(w)}}\Big\|_{L^{\Phi}(w)}{}\nonumber\\{}&&\lesssim \sum_{j=0}^{\infty}\|\langle f,\psi_{2^{j}N}\rangle\psi_{2^{j}N}\|_{L^{\Phi}(w)}h_{\varphi}^{+}(2^{j}N)\end{eqnarray}where in the last inequality we have used Theorem \ref{thmdemFuncOrlpes}. Now using that $\frac{h_{\varphi}^{+}(k)}{k}$ is non-increasing (this follows from  the fact that $\frac{\varphi(t)}{t}$ is  is non-increasing  for all $t> 0,$ see \cite{Benett})  and the definition of the Marcinkiewicz space we have\begin{eqnarray}\label{proofPropJacIneq2} &&\sum_{j=0}^{\infty}\|\langle f,\psi_{Q_{2^{j}N}}\psi_{Q_{2^{j}N}}\|_{L^{\Phi}(w)}h_{\varphi}^{+}(2^{j}N) =\sum_{j=0}^{\infty}\sum_{2^{j-1}N\leq k<2^{j}N}\|\langle f,\psi_{Q_{2^{j}N}}\psi_{Q_{2^{j}N}}\|_{L^{\Phi}(w)}\frac{h_{\varphi}^{+}(2^{j}N)}{2^{j-1}N} {}\nonumber\\ {}&&\leq 2\sum_{k>\frac{N}{2}}\|\langle f,\psi_{Q_{k}}\rangle\psi_{Q_{k}}\|_{L^{\Phi}(w)}\frac{h_{\varphi}^{+}(k)}{k} \leq C\|f\|_{\mathbb{M}_{k^{\alpha}h_{\varphi}^{+}(k)}(\mathcal{B},L^{\Phi}(w))}\sum_{k>\frac{N}{2}}k^{-\alpha}\frac{1}{k} {}\nonumber\\{}&&\leq CN^{-\alpha}\|f\|_{\mathbb{M}_{k^{\alpha}h_{\varphi}^{+}(k)}(\mathcal{B},L^{\Phi}(w))}.\end{eqnarray}\end{proof}
The previous result can be translated as the following inclusion for approximation spaces\begin{eqnarray}\label{icluApprSpaMarc}\mathbb{M}_{k^{\alpha}h_{\varphi}^{+}(k)}(\mathcal{B},L^{\Phi}(w))\hookrightarrow \mathcal{A}_{\infty}^{\alpha}(\mathcal{B},L^{\Phi}(w)).\end{eqnarray}

\subsection{Bernstein type inequalities} Bernstein type estimates are useful to obtain the right hand inclusions for approximation spaces of Theorem \ref{thmInclAprOrlpesoSpa}.
\begin{proposition}\label{PropBernsIneSequSpaces} Let $\Phi$ be a Young function such that $0<i_{\varphi}\leq I_{\varphi}<1,$ $w\in A_{p^{\Phi}}$ a weight in $\mathbb{R}^{d},$ and let $\alpha > 0.$  Let $\mathcal{B}$ be an admissible basis. Then, there exists $C> 0$ so that, for all $N\geq 1,$ and all $f\in\Sigma_{N}$\begin{eqnarray}\label{BernsIneSequSpaces1}\|f\|_{\Lambda_{k^{\alpha}h_{\varphi}^{-}(k)}(\mathcal{B},L^{\Phi}(w))}\leq CN^{\alpha}\|f\|_{L^{\Phi}(w)}.\end{eqnarray}\end{proposition}\begin{proof}\label{proofPropBernIneSequ} Let $f = \sum_{j=1}^{N}\langle f,\psi_{Q_{j}}\psi_{Q_{j}}\in\Sigma_{N},$ written in such a way that $\|\langle f,\psi_{Q_{1}}\psi_{Q_{1}}\|_{L^{\Phi}(w)}\geq \|\langle f,\psi_{Q_{2}}\rangle\psi_{Q_{2}}\|_{L^{\Phi}(w)}\geq\cdots.$  For $1\leq k\leq N,$ using Theorem \ref{democracyinOrlicz1} \begin{eqnarray}\label{proofPropBernIneSequ1}\|\langle f,\psi_{Q_{k}}\psi_{Q_{k}}\|_{L^{\Phi}(w)}h_{\varphi}^{-}(k)&\leq&  C\|\langle f,\psi_{Q_{k}}\rangle \psi_{Q_{k}}\|_{L^{\Phi}(w)}\Big\|\sum_{j=1}^{k}\frac{\psi_{Q_{j}}}{\|\psi_{Q_{j}}\|_{L^{\Phi}(w)}}\Big\|_{L^{\Phi}(w)}{}\nonumber\\{}&&\leq C\|G_{N}(f)\|_{L^{\Phi}(w)}.\end{eqnarray} By (\ref{proofPropBernIneSequ1}) we have \begin{eqnarray*}\|f\|_{\Lambda_{k^{\alpha}h_{\varphi}^{-}(k)}(\mathcal{B},L^{\Phi}(w))} &=& \sum_{k=1}^{N}k^{\alpha}h_{\varphi}^{-}(k)\|\langle f,\psi_{Q_{k}}\rangle \psi_{Q_{k}}\|_{L^{\Phi}(w)}\frac{1}{k}\leq C\|G_{N}(f)\|_{L^{\Phi}(w)}\sum_{k=1}^{N}\frac{k^{\alpha}}{k} \nonumber\\{}&&\leq C'N^{\alpha}\|f\|_{L^{\Phi}(w)}.\end{eqnarray*}\end{proof}
As before, the above result can be stated as an inclusion for approximation spaces. Below, the number $\rho_{\alpha}\in (0,1]$ is chosen so that the quasi-normed space $\Lambda_{k^{\alpha}h_{\varphi}^{-}(k)}$ satisfies the $\rho_{\alpha}$-triangle inequality, that is,  \begin{eqnarray}\label{rhoTrianIne}\|\textbf{s}_{1}+\textbf{s}_{2}\|_{\Lambda_{k^{\alpha}h_{\varphi}^{-}(k)}}^{\rho_{\alpha}}\leq\|\textbf{s}_{1}\|_{\Lambda_{k^{\alpha}h_{\varphi}^{-}(k)}}^{\rho_{\alpha}}+\|\textbf{s}_{2}\|_{\Lambda_{k^{\alpha}h_{\varphi}^{-}(k)}}^{\rho_{\alpha}}.\end{eqnarray}
\begin{corollary}\label{CoricluApprSpaLorSpa}Let $\alpha > 0.$ Then, with the same hypothesis as in Proposition \ref{PropBernsIneSequSpaces}, we have \begin{eqnarray}\label{CoricluApprSpaLorSpa1}\mathcal{A}_{\rho_{\alpha}}^{\alpha}(\mathcal{B},L^{\Phi}(w))\hookrightarrow \Lambda_{k^{\alpha}h_{\varphi}^{-}(k)}(\mathcal{B},L^{\Phi}(w))\end{eqnarray} \end{corollary}\begin{proof}\label{proofCorincluApprSpaLor} The argument for (\ref{CoricluApprSpaLorSpa1}) is standard (see, e.g., \cite{DP}). It suffices to prove that $$\|f\|_{\Lambda_{k^{\alpha}h_{\varphi}^{-}(k)}(\mathcal{B},L^{\Phi}(w))}\leq C\|f\|_{\mathcal{A}_{\rho_{\alpha}}^{\alpha}(\mathcal{B},L^{\Phi}(w))}, \quad \forall f\in\Sigma_{N},\ N\geq 1$$ with a constant $C> 0$ independent of $N$ and one obtains the desired inclusion by letting $N\longrightarrow \infty.$ We also assume that $N = 2^{J}.$ Now, write $f= \sum_{j=0}^{J}[f^{(j)}-f^{(j-1})],$ where by convection $f^{(J)} = f, \ f^{(-1)} = 0$ and $f^{(j)}\in \Sigma_{2^{j}}$ is so that $\|f - f^{(j)}\|_{L^{\Phi}(w)}\leq 2\sigma_{2^{j}}(f)_{L^{\Phi}(w)}, \ 0\leq j<J.$ Then applying, (\ref{rhoTrianIne}) and  Proposition \ref{PropBernsIneSequSpaces} to $f^{(j)}-f^{(j-1)}\in \Sigma_{2^{j+1}}$ we obtain \begin{eqnarray*}\|f\|_{\Lambda_{k^{\alpha}h_{\varphi}^{-}(k)(\mathcal{B}, L^{\Phi}(w))}}&\leq& \Big[\sum_{j=0}^{J}\|f^{(j)}-f^{(j-1)}\|_{\Lambda_{k^{\alpha}h_{\varphi}^{-}(k)}(\mathcal{B}, L^{\Phi}(w))}^{\rho}\Big]^{\frac{1}{\rho}} {}\nonumber\\{}&&\leq C\Big[\sum_{j=0}^{J}2^{j\alpha\rho}\|f^{(j)}-f^{(j-1)}\|_{L^{\Phi}(w)}^{\rho}\Big]^{\frac{1}{\rho}}.\end{eqnarray*}
Now, by assumption, for $1\leq j\leq J$ $$\|f^{(j)}-f^{(j-1)}\|_{L^{\Phi}(w)}\leq \|f^{(j)}-f\|_{L^{\Phi}(w)} + \|f - f^{(j-1)}\|_{L^{\Phi}(w)}\leq 4\sigma_{2^{(j-1)}}(f)_{L^{\Phi}(w)}.$$ On the other hand, for $j=0$ we have $$\|f^{(0)}-f^{(-1)}\|_{L^{\Phi}(w)} = \|f^{(0)}\|_{L^{\Phi}(w)}\leq\|f^{(0)}-f\|_{L^{\Phi}(w)}+\|f\|_{L^{\Phi}(w)}\leq2\sigma_{1}(f)_{L^{\Phi}(w)}+\|f\|_{L^{\Phi}(w)}.$$
Hence,

\begin{eqnarray*} \|f\|_{\Lambda_{k^{\alpha}h_{\varphi}^{-}(k)}(\mathcal{B},L^{\Phi}(w))}\leq C\Big[\|f\|_{L^{\Phi}(w)}+\sum_{j=0}^{J-1}(2^{j\alpha}\sigma_{2^{j}}(f)_{L^{\Phi}(w)})^{\rho}\Big]^{\frac{1}{\rho}}\approx \|f\|_{\mathcal{A}_{\rho}^{\alpha}(\mathcal{B},L^{\Phi}(w))}.\end{eqnarray*}\end{proof}

Finally, using real interpolation we can obtain inclusions for the whole family of approximation spaces $\mathcal{A}_{q}^{\alpha}(\mathcal{B},L^{\Phi}(w)), \ 0<q\leq\infty.$ For this we consider the interpolation properties of the sequence spaces $\Lambda_{\eta}^{q},$ namely,\begin{eqnarray}\label{interpLorenSpaces}(\Lambda_{k^{\alpha_{0}}\eta(k)}^{r},\Lambda_{k^{\alpha_{1}}\eta(k)}^{r})_{\alpha,q}=\Lambda_{k^{\alpha}\eta(k)}^{q},\quad \alpha = (1-\theta)\alpha_{0}+\theta\alpha_{1},\end{eqnarray} for all $0<q,r\leq\infty,\ 0<\theta<1$ (see, e.g., \cite{Person} Proposition 6.2, \cite{Meruci}, Theorem 3).
\begin{theorem}\label{ThmInclWeiDiscOrliAppSpac}Let $\Phi$ be a Young function such that $0<i_{\varphi}\leq I_{\varphi}<1,$ $w\in A_{p^{\Phi}}$ a weight in $\mathbb{R}^{d},$  $\alpha > 0,$ and $ 0<q\leq\infty.$ Let $\mathcal{B}$ be an admissible wavelet basis. Then \begin{eqnarray}\label{ThmInclWeiOrliAppSpac1}\Lambda_{k^{\alpha}h_{\varphi}^{+}(k)}^{q}(\mathcal{B},L^{\Phi}(w))\hookrightarrow\mathcal{A}_{q}^{\alpha}(\mathcal{B},L^{\Phi}(w))\hookrightarrow\Lambda_{k^{\alpha}h_{\varphi}^{-}(k)}^{q}(\mathcal{B},L^{\Phi}(w)).\end{eqnarray}
\end{theorem}\begin{proof}Let $\alpha_{0}<\alpha<\alpha_{1},$ so that $\alpha = (\alpha_{0} + \alpha_{1})/2.$ Then, for every $0<q,r\leq\infty$ we have (see, e.g., \cite{DP}) $$ \mathcal{A}_{q}^{\alpha}(\mathcal{B},L^{\Phi}(w)) = (\mathcal{A}_{r}^{\alpha_{0}}(\mathcal{B},L^{\Phi}(w)),\mathcal{A}_{r}^{\alpha_{1}}(\mathcal{B},L^{\Phi}(w)))_{\frac{1}{2},q}.$$ Letting $r = \min(\rho_{\alpha_{0}},\rho_{\alpha_{1}})$ and using (\ref{CoricluApprSpaLorSpa1})\begin{eqnarray*} \mathcal{A}_{q}^{\alpha}(\mathcal{B},L^{\Phi}(w))&=&(\mathcal{A}_{r}^{\alpha_{0}}(\mathcal{B},L^{\Phi}(w)),\mathcal{A}_{r}^{\alpha_{1}}(\mathcal{B},L^{\Phi}(w)))_{\frac{1}{2},q}{}\nonumber\\{}&&\hookrightarrow (\Lambda_{k^{\alpha_{0}}h_{\varphi}^{-}(k)}(\mathcal{B},L^{\Phi}(w)),\Lambda_{k^{\alpha_{1}}h_{\varphi}^{-}(k)}(\mathcal{B},L^{\Phi}(w)))_{\frac{1}{2},q} {}\nonumber\\{}&&= \Lambda_{k^{\alpha}h_{\varphi}^{-}(k)}^{q}(\mathcal{B},L^{\Phi}(w)),\end{eqnarray*} where the last equality follows from (\ref{interpLorenSpaces}). Similarly, by (\ref{icluApprSpaMarc})\begin{eqnarray*}\mathcal{A}_{q}^{\alpha}(\mathcal{B},L^{\Phi}(w)) &=& (\mathcal{A}_{\infty}^{\alpha_{0}}(\mathcal{B},L^{\Phi}(w)),\mathcal{A}_{\infty}^{\alpha_{1}}(\mathcal{B},L^{\Phi}(w)))_{\frac{1}{2},q}{}\nonumber\\{}&&\hookleftarrow (\mathbb{M}_{k^{\alpha_{0}}h_{\varphi}^{+}(k)}(\mathcal{B},L^{\Phi}(w)),\mathbb{M}_{k^{\alpha_{1}}h_{\varphi}^{+}(k)}(\mathcal{B},L^{\Phi}(w))){}\nonumber\\{}&&=\Lambda_{k^{\alpha}h_{\varphi}^{+}(k)}^{q}(\mathcal{B},L^{\Phi}(w)).\end{eqnarray*}   \end{proof}

We now prove that the inclusions (\ref{ThmInclWeiOrliAppSpac1}) are optimal. To state the Theorem we write $\mathrm{D}$ for the class of sequences $\eta = \{\eta(k)\}_{k=1}^{\infty}$  that are increasing,   and  doubling.
\begin{theorem}\label{thmOptimIncl} Same hypothesis as in Theorem \ref{ThmInclWeiDiscOrliAppSpac}. For fixed $\alpha>0$ and $q,\ 0<q\leq\infty,$ the inclusions given in (\ref{ThmInclWeiOrliAppSpac1}) are best possible in the scale of weighted Lorentz spaces $\Lambda_{k^{\alpha}\eta(k)}^{q}(\mathcal{B},L^{\Phi}(w))$ where $\eta\in\mathrm{D}$\end{theorem}\begin{proof}Suppose $\Lambda_{k^{\alpha}\eta(k)}^{q}(\mathcal{B},L^{\Phi}(w))\hookrightarrow \mathcal{A}_{q}^{\alpha}(\mathcal{B},L^{\Phi}(w)).$ We want to prove that $h_{\varphi}^{+}(N)\lesssim \eta(N)$ for all $N=1,2,\ldots.$ By definition of $h_{\varphi}^{+}(N)$ we can choose $\tau = \tau(N)>0$ such that

\begin{eqnarray}\label{proofThmOptiIncl} \frac{\varphi(N\tau)}{\varphi(\tau)}\leq h_{\varphi}^{+}(N)\leq 2\frac{\varphi(N\tau)}{\varphi(\tau)}.\end{eqnarray} 

By Lemma \ref{auxLemm2} we can choose a sequence of pairwise disjoint cubes $\Gamma = \{R_{j}\}_{j=1}^{2N}$ such that $w(R_{j})\approx \tau.$ Let $\widetilde{1}_{\Gamma} = \sum_{j=1}^{2N}\frac{\psi_{R_{j}}}{\|\psi_{R_{j}}\|_{L^{\Phi}(w)}}.$ By Theorem \ref{thmcaraOrlPes}, $\|f\|_{L^{\Phi}(w)}$ is equivalent to the lattice norm $\|S_{\psi}(f)\|_{L^{\Phi}(w)};$ thus there exists $\Gamma'\subset\Gamma$ with $\Gamma'= N$ such that $\sigma_{N}(\widetilde{1}_{\Gamma})\approx \|\widetilde{1}_{\Gamma'}\|_{L^{\Phi}(w)}$ (see (2.6) in \cite{GustavoSharpJackson}). Thus, by (\ref{proof-PropDicCub}) and (\ref{proofThmOptiIncl})$$\sigma_{N}(\widetilde{1}_{\Gamma})_{L^{\Phi}(w)}\approx\|\widetilde{1}_{\Gamma'}\|_{L^{\Phi}(w)}\approx\frac{\varphi(N\tau)}{\varphi(\tau)}\approx h_{\varphi}^{+}(N).$$ Hence,

\begin{eqnarray}\label{proofThmOptiIncl1}\|\widetilde{1}_{\Gamma}\|_{\mathcal{A}_{q}^{\alpha}(\mathcal{B},L^{\Phi}(w))}\geq\Big(\sum_{k=N/2}^{N}k^{\alpha q}\sigma_{k}(\widetilde{1}_{\Gamma})^{q}\frac{1}{k}\Big)^{\frac{1}{q}}\gtrsim \sigma_{N}(\widetilde{1}_{\Gamma})N^{\alpha}\approx N^{\alpha}h_{\varphi}^{+}(N).\end{eqnarray}  On the other hand

\begin{eqnarray}\label{proofThmOptiIncl2} \|\widetilde{1}_{\Gamma}\|_{\Lambda_{k^{\alpha}\eta(k)}^{q}(\mathcal{B},L^{\Phi}(w))} = \Big(\sum_{k=1}^{2N}(k^{\alpha}\eta(k))^{q}\frac{1}{k}\Big)^{\frac{1}{q}}\lesssim\eta(2N)N^{\alpha}\lesssim \eta(N)N^{\alpha}\end{eqnarray} by the doubling property of $\eta.$ The inequalities (\ref{proofThmOptiIncl1}) and (\ref{proofThmOptiIncl2}) together with  our assumption imply the desired result.\\

Suppose now that $\mathcal{A}_{q}^{\alpha}(\mathcal{B},L^{\Phi}(w))\hookrightarrow \Lambda_{k^{\alpha}\eta(k)}^{q}(\mathcal{B},L^{\Phi}(w)).$ We want to prove that $\eta(k)\leq h_{\varphi}^{-}(N)$ for all $N=1,2,\ldots.$ Let $\Gamma\subset\mathcal{D}$ with $|\Gamma|= N.$ Write $\widetilde{1}_{\Gamma} = \sum_{Q\in\Gamma}\frac{\psi_{Q}}{\|\psi_{Q}\|_{L^{\Phi}(w)}}.$ Since $\sigma_{k}(\widetilde{1}_{\Gamma}) \leq \|\widetilde{1}_{\Gamma}\|_{L^{\Phi}(w)}$ for all $k=1,2,\ldots,N,$ our hypothesis imply 

\begin{eqnarray}\label{proofThmOptiIncl3} \|\widetilde{1}_{\Gamma}\|_{\Lambda_{k^{\alpha}\eta(k)}^{q}(\mathcal{B},L^{\Phi}(w))}\lesssim \|\widetilde{1}_{\Gamma}\|_{\mathcal{A}_{q}^{\alpha}(\mathcal{B},L^{\Phi}(w))}\lesssim N^{\alpha}\|\widetilde{1}_{\Gamma}\|_{L^{\Phi}(w)}.\end{eqnarray} 

On the other hand 

\begin{eqnarray}\label{proofThmOptiIncl4} \|\widetilde{1}_{\Gamma}\|_{\Lambda_{k^{\alpha}\eta(k)}^{q}(\mathcal{B},L^{\Phi}(w))}\geq \Big(\sum_{k=\frac{N}{2}}^{N}(\eta(k)k^{\alpha})^{q}\frac{1}{k}\Big)^{\frac{1}{q}}\gtrsim N^{\alpha}\eta(N/2)\gtrsim N^{\alpha}\eta(N)\end{eqnarray} 

since $\eta$ is doubling. By (\ref{proofThmOptiIncl3}) and (\ref{proofThmOptiIncl4}) we have $\eta(N)\lesssim \|\widetilde{1}_{\Gamma}\|_{L^{\Phi}(w)}$ for all $\Gamma\subset\mathcal{D},$ with $|\Gamma| = N.$ Taking the infimum  over all $\Gamma\in\mathcal{D},$ with $|\Gamma| = N,$ we obtain $\eta(N)\leq h_{l}(N;\mathcal{B},L^{\Phi}(w))\approx h_{\varphi}^{-}(N)$ by Theorem \ref{thmdemFuncOrlpes}.
  \end{proof}

 \section{Approximation spaces for $L^{p}(w)$}\label{secAppSpcLppes}

Corollary \ref{coro-CarLpesAppr} is now an easy consequence of Theorem \ref{ThmInclWeiDiscOrliAppSpac} and Corollary \ref{CorodemocFuncLppeso}. The rest of this section is devoted to prove Theorem \ref{thmCaraLppeso} ( see Theorem \ref{thmcaractLppesApprSpac} bellow).

The  approximation spaces $\mathcal{A}_{\tau}^{\gamma}(\mathcal{B},L^{p}(w))$ can also be identified with weighted  Besov  spaces. Our definition of weighted
 Besov spaces is borrowed from \cite{Rou1,Rou2}, and it is
 modeled on the corresponding definition of Besov spaces without
 weights developed in \cite{Pee} (see also \cite{FJ1} and
 \cite{FJ2}).\\

 We say that a function $\varphi \in \mathscr{S}(\R^{d})$ belongs
 to the class of admissible kernels if $Supp \widehat{\varphi}\subset \{\xi \in \R^{d}:
 \frac{1}{2}<|\xi|<2\}$ and $|\widehat{\varphi}(\xi)|\geq c > 0,$  if
 $\frac{3}{5}<|\xi|<\frac{5}{3}.$ Set $\varphi_{k}(x) =
 2^{kd}\varphi(2^{k}x)$ for $k \in \mathbb{Z}.$\\

 Let $\alpha \in\mathbb{R},\ 1\leq p < \infty,\  0< q \leq \infty,$
 $\varphi$ admissible kernel, and $w$ an $A_{p}$ weight on $\R^{d}.$ The
 homogeneous weighted Besov space $\dot{B}_{p,q}^{\alpha}(w)$ is the
 set of all tempered distributions $f\in \mathscr{S}'/\mathscr{P}$
 (modulo polynomials) such that
 \begin{eqnarray}\label{norm-WeightBesovSpace}\|f\|_{\dot{B}_{p,q}^{\alpha}(w)}
 =
 \Big[\sum_{k\in\mathbb{Z}}(2^{k\alpha}\|\varphi_{k}*f\|_{L^{p}(w)})^{q}\Big]^{1/q}<\infty.\end{eqnarray}

 This definition depends initially of the choice of admissible
 $\varphi.$ It can be proved (see Theorem 1.8 in \cite{Rou1} or
 \cite{Rou2}) that this is independent of the choice of admissible
 $\varphi.$ Also, the spaces $\dot{B}_{p,q}^{\alpha}(w)$ are
 (quasi)-Banach spaces (see section 4.4 of \cite{Rou2}). \\

 Let $\Psi = \{\psi^{l}: l=1,2,\ldots,2^{d-1}\}$ be an orthonormal wavelet family in $L^{2}(\mathbb{R}^{d})$ constructed from the 1-dimensional Lemari\'e-Meyer wavelets (see \cite{Weiss,LemarieMeyer, Me}). Write $s_{Q}^{l} = \langle f,\psi_{Q}^{l}\rangle,\ Q\in \mathcal{D},\ l=1,2,\ldots,2^{d-1}$ for the wavelet coefficients.

\begin{proposition}\label{carac-WeigBesovSpace}(see Theorem 10.2 in \cite{Rou1} or Theorem 6.2 in
 \cite{Rou2}).

  Let $\alpha \in\R, \ 0<q\leq \infty, \ 1\leq p< \infty$ and let $w$ be an $A_{p}$ weight in $\mathbb{R}^{d}.$  Let $\Psi$ be a family of
  Lemari\'e-Meyer wavelets as defined above.
  Then\begin{eqnarray}\label{carac-WeigBesovSpace1}\|f\|_{\dot{B}_{p,q}^{\alpha}(w)}&\approx&
  \sum_{l=1}^{2^{d}-1}\Big[\sum_{j\in\mathbb{Z}}\Big(\sum_{|Q|=  2^{-jd}}(|Q|^{-\frac{\alpha}{d}-\frac{1}{2}}|s_{Q}^{l}|w(Q)^{\frac{1}{\tau}})^{\tau}\Big)^{\frac{q}{\tau}}\Big]^{\frac{1}{q}}.\end{eqnarray}
 \end{proposition}\begin{Remark}\label{rem-carac-WeigBesovSpac} It is also proved in Theorem 10.2 of \cite{Rou1} and Theorem 6.2 of \cite{Rou2}
 that the condition $w$ doubling, that is, there exists $C> 0$ such that $$\int_{B_{2\delta}(z)}w(x)dx\leq C\int_{B_{\delta}(z)}w(x)dx,\quad \forall z\in\R^{d}\ and\ \forall \delta > 0,$$
 is sufficient to guarantee the equivalence (\ref{carac-WeigBesovSpace1}).
 \end{Remark} 
 
 \begin{Remark}\label{rem-carac-WeigBesovSpac1}
 Equivalence (\ref{carac-WeigBesovSpace1}) also holds for the
 family \
$_{\bf{N}}\Psi = \{_{N}\psi^{l}: l=1,\ldots,2^{d-1}\}$ constructed
from the 1-dimensional Daubechies compactly supported wavelets (see
\cite{Dau}), provided $N$ is sufficiently large (see \cite{Rou1}).
\end{Remark}
 \begin{theorem}\label{thmcaractLppesApprSpac} Let $\gamma >
 0$,  $1 <  p< \infty.$  We have
 \begin{eqnarray}\label{caractEsapacioAproximacion}
  \mathcal{A}_{\tau}^{\gamma/d}(\Psi, L^{p}(w)) =l^{\tau}(\Psi,L^{p}(w)) = \dot{B}_{\tau,\tau}^{\gamma}(w^{\frac{\tau}{p}}),
 \quad whenever \quad \frac{1}{\tau} =
 \frac{\gamma}{d} + \frac{1}{p},\end{eqnarray} for all $w\in A_{\tau}(\mathbb{R}^{d})$ and all orthonormal wavelet families $\Psi$ for which (\ref{carac-WeigBesovSpace1}) holds for $\dot{B}_{\tau,\tau}^{\gamma/d}(w).$ \end{theorem}
 For the proof we shall need the following lemma:
 \begin{lemma}\label{equivalenciapesos} Let $w \in A_{r}$ be a weight in
 $\R^d$, $r \geq  1$,  $0< \delta  < 1,$  and $u(x) = w(x)^{\delta}$.
 Then   $u \in A_{r}$, and $w_{Q}\approx (u_{Q})^{\frac{1}{\delta}}$,
 where \begin{eqnarray*}w_{Q}=\frac{1}{|Q|}\int_{Q}w(x)dx.\end{eqnarray*}\end{lemma}
  \begin{proof} When  $r>1$, since $w \in A_{r}$ and $0<\delta<1$, using Jensen's
  inequality we have\setlength\arraycolsep{2pt}
  \begin{eqnarray*}&&\Big(\frac{1}{|Q|}\int_{Q}u(x)dx\Big)\Big(\frac{1}{|Q|}\int_{Q}u^{1-r'}(x)dx\Big)^{r-1}\nonumber
  \\
  &=& \Big(\frac{1}{|Q|}\int_{Q} w(x)^{\delta} dx \Big)\Big(\frac{1}{|Q|}\int_{Q}w(x)^{\delta(1-r')}dx\Big)^{r-1}{}
  \nonumber\\ {}
  &\leq& \Big[\Big(\frac{1}{|Q|}\int_{Q}w(x)dx\Big)\Big(\frac{1}{|Q|}\int_{Q}w(x)^{1-r'}dx\Big)^{r-1}\Big]^{\delta} \leq C_{w}^{\delta}.\end{eqnarray*} Thus, $u\in A_{r}$.
  Now we shall to prove the equivalence $w_{Q}\approx(u_{Q})^{\frac{1}{\delta}}$.
   On the one hand, using  Jensen's inequality with $\delta < 1$, we have
   \begin{eqnarray*}(u_{Q})^{\frac{1}{\delta}}=\Big(\frac{1}{|Q|}
   \int_{Q}w(x)^{\delta}dx\Big)^{\frac{1}{\delta}}\leq\Big(\frac{1}{|Q|}\int_{Q}w(x)dx\Big) = w_{Q}.\end{eqnarray*}
    On the other hand, since  $h(t) = t^{-(r-1)\delta},\ t > 0,$ is a convex function, using again
    Jensen's inequality we have
    \begin{eqnarray*}\Big(\frac{1}{|Q|}\int_{Q}w^{1-r'}(x)dx\Big)^{-(r-1)\delta}
    \leq\Big(\frac{1}{|Q|}\int_{Q}w(x)^{\delta}dx\Big)= u_{Q}.\end{eqnarray*}
  From the condition $w\in A_{r} $  it follows that
  \begin{eqnarray*}w_{Q}=\Big(\frac{1}{|Q|}\int_{Q}w(x)dx\Big)\leq C_{w} \Big(\frac{1}{|Q|}\int_{Q}w^{1-r'}(x)dx\Big)^{-(r-1)}\leq C(u_{Q})^{\frac{1}{\delta}}.  \end{eqnarray*}
  For $r = 1$, on the one hand, since $w \in A_{1}$, for almost all $x \in Q$, using  again  Jensen's inequality
   we have
   \begin{eqnarray*}\Big(\frac{1}{|Q|}\int_{Q}u(x)dx\Big) =
   \Big(\frac{1}{|Q|}\int_{Q}w^{\delta}(x)dx\Big)\leq\Big(\frac{1}{|Q|}\int_{Q}w(x)dx\Big)^{\delta}
   \leq C w(x)^{\delta}= C u(x).\end{eqnarray*}
  Thus $u\in A_{1}$. On the other hand  we can use again    Jensen's inequality and obtain, $(u_{Q})^{\frac{1}{\delta}}\leq w_{Q}.$ Moreover, the condition $w \in A_{1}$  implies that
  \begin{eqnarray*}w_{Q}=\frac{1}{|Q|}\int_{Q}w(x)dx \leq C ess_{Q}\inf w = C(ess_{Q}\inf u)^{\frac{1}{\delta}}\leq C\Big(\frac{1}{|Q|}\int_{Q}u(x)dx\Big)^{\frac{1}{\delta}}=
  C(u_{Q})^{\frac{1}{\delta}}\end{eqnarray*}\end{proof}
  \begin{proof}(\textbf{of Theorem
  \ref{thmcaractLppesApprSpac}}).  The  first equality   in  (\ref{caractEsapacioAproximacion}) follows  from Corollary \ref{coro-CarLpesAppr} (with $\tau = q$ and $\alpha = \frac{\gamma}{d}$). For the second equality,
  observe that for a single element of the basis $\Psi,$  we have \begin{eqnarray}\label{norLppesoSingEl}\|\psi_{Q}\|_{L^{p}(w)}&\approx& \Big\|\Big(\frac{\chi_{Q}(x)}{|Q|}\Big)^{\frac{1}{2}}\Big\|_{L^{p}(w)} = |Q|^{-\frac{1}{2}}\|\chi_{Q}(x)\|_{L^{p}(w)} = |Q |^{-\frac{1}{2}}\Big(\int_{Q}w(x)dx\Big)^{\frac{1}{p}}{}\nonumber\\{}&&=|Q|^{-\frac{1}{2}}w(Q)^{\frac{1}{p}}.\end{eqnarray}

Let $u(x) = w(x)^{\frac{\tau}{p}}.$ By  Lemma
\ref{equivalenciapesos} with $r =\tau$ and $\delta = \frac{\tau}{p} <
1$ we deduce that $u\in A_{\tau}\subset A_{p}$  and $w_{Q}^{\frac{1}{p}} \approx
(u_{Q})^{\frac{1}{\tau}}.$ Thus, since $w(Q) = |Q|w_{Q}$ we obtain
\begin{eqnarray*}
\|f\|_{\ell^{\tau}(\Psi,L^{p}(w))} &=& \Big\|\|\langle f,\psi_{Q_{k}}\rangle\psi_{Q_{k}}\|_{L^{p}(w)}\Big\|_{\ell^{\tau}} \approx
\Big(\sum_{Q\in\mathcal{D}}(|Q|^{-\frac{1}{2}+\frac{1}{p}}|\langle f,\psi_{Q}\rangle|(w_{Q})^{\frac{1}{p}})^{\tau}\Big)^{\frac{1}{\tau}}
{}\nonumber\\{}&&\approx
\Big(\sum_{Q\in\mathcal{D}}(|Q|^{-\frac{1}{2}+\frac{1}{p}}|\langle f,\psi_{Q}\rangle|(u_{Q})^{\frac{1}{\tau}})^{\tau}\Big)^{\frac{1}{\tau}}
 {}\nonumber\\{}&&=\Big(\sum_{Q\in\mathcal{D}}(|Q|^{-\frac{1}{2}-\frac{\gamma}{d} +
 \frac{1}{\tau}}|\langle f,\psi_{Q}\rangle|(u_{Q})^{\frac{1}{\tau}})^{\tau}\Big)^{\frac{1}{\tau}}{} {}\nonumber\\{}&&=
 \Big(\sum_{Q\in\mathcal{D}}(|Q|^{-\frac{1}{2}-\frac{\gamma}{d}}
 |\langle f,\psi_{Q}\rangle|(u(Q))^{\frac{1}{\tau}})^{\tau}\Big)^{\frac{1}{\tau}}
 \approx\|f\|_{\dot{B}_{\tau,\tau}^{\gamma}(w^{\frac{\tau}{p}})}.
 \end{eqnarray*} \end{proof}As a corollary we prove a non-trivial
 interpolation result.\begin{corollary}\label{corol-Interpol-AppS} Let $\gamma>0,\ 1<p<\infty,$   $\frac{1}{\tau} = \frac{\gamma}{d} + \frac{1}{p},$ and $w \in A_{\tau}(\mathbb{R}^{d}).$ Let $\Psi$ be an orthonormal wavelet family for which (\ref{carac-WeigBesovSpace1}) holds for the Besov spaces involved in this Corollary.   For  $0<\theta <1$ we have
 $$ (L^{p}(w),\dot{B}_{\tau,\tau}^{\gamma}(w^{\tau/p}))_{\theta,\tau_{\theta}} = \dot{B}_{\tau_{\theta},\tau_{\theta}}^{\theta\gamma}(w^{\tau_{\theta}/p})$$
 where  $\frac{1}{\tau_{\theta}} = \frac{\theta\gamma}{d} + \frac{1}{p}.$\end{corollary}
 \begin{proof} If $\Phi(t)=t^{p},$  use Proposition \ref{PropJacInqSeqOrlSpa},  the
 continuous embedding $l^{\tau}\hookrightarrow l^{\tau,\infty},$ and
 Theorem \ref{thmcaractLppesApprSpac} to obtain, for all
 $N=1,2\ldots$
 \setlength\arraycolsep{2pt}\begin{eqnarray}\label{proof-Corol-InterpApS1}\sigma_{N}(
 f)_{L^{p}(w)}&\leq& CN^{-(\frac{1}{\tau}-\frac{1}{p})}\| f\|_{l^{\tau,\infty}(\Psi,L^{p}(w))}\leq C N^{-(\frac{1}{\tau}-\frac{1}{p})}\|f\|_{l^{\tau}(\Psi,L^{p}(w))}{}
 \nonumber\\ {} \leq CN^{-(\frac{1}{\tau}-\frac{1}{p})}\|
 f\|_{\dot{B}_{\tau,\tau}^{\gamma}(w^{\tau/p})},\end{eqnarray} whenever $\frac{1}{\tau} = \frac{\gamma}{d} + \frac{1}{p}.$

 From  Theorem \ref{thmcaractLppesApprSpac}, and Proposition \ref{PropBernsIneSequSpaces}
  we obtain, for all $ g \in \Sigma_{N},\
 N=1,2,3,\ldots$ \begin{eqnarray}\label{proof-Corol-InterpApS2} \|
 g\|_{\dot{B}_{\tau,\tau}^{\gamma}(w^{\tau/p})}\leq C\|
 g\|_{l^{\tau}(\Psi,L^{p}(w))}\leq CN^{\frac{1}{\tau}-\frac{1}{p}}\|
 g\|_{L^{p}(w)},\end{eqnarray} whenever $\frac{1}{\tau} = \frac{\gamma}{d} +
 \frac{1}{p}.$ From (\ref{proof-Corol-InterpApS1}),
 (\ref{proof-Corol-InterpApS2}), and the general theory developed by R.
 DeVore and V. A. Popov (see Theorem 3.1 in \cite{DP}) we deduce
 \begin{eqnarray}\label{proof-Corol-InterpApS3} \mathcal{A}_{q}^{\frac{\theta\gamma}{d}}(\Psi,L^{p}(w)) = (L^{p}(w),\dot{B}_{\tau,\tau}^{\gamma}(w^{\tau/p}))_{\theta,q} \end{eqnarray}
 whenever $0<q\leq \infty$ and $0<\theta < 1.$ We use again Theorem \ref{thmcaractLppesApprSpac}
 to obtain $$ \mathcal{A}_{\tau_{\theta}}^{\frac{\theta\gamma}{d}}(\Psi,L^{p}(w)) = \dot{B}_{\tau_{\theta},\tau_{\theta}}^{\theta\gamma}(w^{\tau_{\theta}/p})$$
 when $\frac{1}{\tau_{\theta}} = \frac{\theta\gamma}{d} + \frac{1}{p}.$ The result follows from (\ref{proof-Corol-InterpApS3}) with $q =\tau_{\theta}.$\end{proof}


\begin{thebibliography}{1}
\bibitem{Aimar}
\textsc{H.A. Aimar, A.L. Bernardis, and F.J. Mart\'in-Reys},
\emph{Multiresolution Approximation and Wavelet Bases of Weighted
Lebesgue Spaces},  J. Fourier Anal. and Appl., 9, nº 5, (2003),
497-510.


\bibitem{Benett}
 \textsc{C. Benett, and R.C. Sharpley},
 \emph{Interpolation of Operators},
  Pure and Appl. Math. 129, Academic Press, 1988.

\bibitem{CRSoria}
\textsc{M.J. Carro, J. Raposo, and J. Soria},
\emph{Recent developments in the theory of Lorentz spaces and weighted inequalities}, Mem. Amer. Math. Soc. 107 (2007).

\bibitem{MatellUribePerez}
\textsc{D. Cruz-Uribe, J.M. Martell, and C. P\'erez},
\emph{Weights, Extrapolation and Theory of Rubio de Francia}, Preprint, (2008).
\bibitem{Dau}
\textsc{I. Daubechies}, \emph{Orthonormal bases of compactly
supported wavelets}, Comm. Pure Appl. Math., 2, (1986), 1--18.

\bibitem{DPetrova}
\textsc{R. DeVore, G. Petrova, and V. Temlyakov},
\emph{Best basis selection for approximation in $L^{p}$},
Found. Comput. Math. 3, (2003), 161--185.
\bibitem{DP}
     \textsc{R. DeVore and   V.A. Popov},
     \emph{Interpolation spaces and nonlinear approximation},
     Function spaces and applications (Lund, 1986),
     Lecture Notes in Math., 1302, Springer, Berlin, (1988),
      191--205.



\bibitem{DKKT}
\textsc{S.J. Dilworth, N.J. Kalton, D. Kutzarova, and V.N. Temlyakov},
\emph{The Thresholding Greedy Algorithm, Greedy Bases, and Duality}, Construtive Approximation, 19, (2003), 575--597.

\bibitem{FJ1}
\textsc{M. Frazier and B. Jawerth}, \emph{Decomposition of Besov
Spaces}, Indiana Univ. Math. J., 34, (1985), 777--799.

\bibitem{FJ2}
\textsc{M. Frazier and B. Jawerth}, \emph{A discrete transform and
decomposition of distribution spaces}, J. Func. Anal., 93, (1990),
34--170.

\bibitem{Rubio}
\textsc{J. Garcia-Cuerva and J.L. Rubio de Francia}, \emph{ Weighted
Norm Inequalities and Related Topics},  North Holland Mathematical
Studies 116, Elsevier Sciences Publishers, 1985.

\bibitem{CUERVAMARTELL}
\textsc{J. Garc\'ia-Cuerva and J.M. Martell},
\emph{Wavelet characterization of weighted spaces},
J. Geom. Anal. 11, no. 2,(2001), 241--262.
\bibitem{GustavoSharpJackson}
\textsc{G. Garrig\'os and E. Hern\'andez}, \emph{ Sharp Jackson and
Bernstein Inequalities for n-term Approximation
  in Sequence Spaces with Applications},
Indiana Univ. Math. J., 53 (2004), 1739--1762.


\bibitem{GHM}
\textsc{G. Garrig\'os, E.Hern\'andez, and J.M. Martell},\emph{
Wavelets, Orlicz and Greedy bases},  Appl. Compt. Harmon. Anal.,
24, (2008), 70--93.


\bibitem{GNielsen}
\textsc{R. Gribonval, and M. Nielsen},
\emph{Some Remarks on non-linear approximation with Schauder bases},
East J. Approx. 7, (3), (2001), 267--285.


\bibitem{Weiss}
\textsc{E. Hern\'andez and G.Weiss}, \emph{ A first course on
wavelets}, CRC Press, Boca Raton FL, 1996.
\bibitem{HJLY}
 \textsc{C. Hsiao, B. Jawerth, B.J. Lucier,  and X.M. Yu},
     \emph{Near optimal compression of almost optimal wavelet
     expansions},
      Wavelets: mathematics and applications,
     Stud. Adv. Math., CRC, Boca Raton, FL, (1994), 425--446.
\bibitem{KT}
\textsc{A. Kamont and V. N. Temlyakov},
\emph{Greedy approximation and the multivariate Haar system}, Studia Math. 161, no. 3,(2004),  199--223.
\bibitem{KerkiacharyanNon}
\textsc{G. Kerkiacharyan, and D. Picard},
\emph{Nonlinear Approximation and Mckenhoupt Weights},
Constr. Approx. 20, (2006), 123--156.

\bibitem{Krebec}
\textsc{V. Kokilashvili, and M. Krebec},
\emph{Weighted Inequalities in Lorentz and Orlicz spaces},
Word Scientific, Singapore, 1991.
\bibitem{Konyagin}
\textsc{S.V. Konyagin and V.N. Temlyakov},
      \emph {A remark on greedy approximation in Banach spaces.}
      East J. Approx., 5, (1999), 365--379.
\bibitem{KPS}
\textsc{S. Krein, J. Petunin, and E. Semenov},
\emph{Interpolation of Linear Operators}, Translations Math. Monographs, vol. 55, Amer. Math. Soc., Providence, RI, 1982.

\bibitem{LemarieMeyer}
\textsc{P.G. Lemari\'e and Y. Meyer}, \emph{Ondelettes et bases
hilbertiannes},  Rev. Mat. Iberoamericana, 2, no. 1-2, (1986), 1--18.
\bibitem{Meruci}
\textsc{C. Meruci},
\emph{Applications of interpolation with a function parameter to Lorentz, Sobolev and Besov spaces}, in : M. Cwikel, J. Peetre (Eds.), \emph{Interpolation spaces and Allied Topics in Analysis}, Lund, 1983, in: Lecture Notes in Math., vol.1070, Springer, Berlin, (1984), pp. 183--201-

\bibitem{Me}
 \textsc{Y. Meyer},
     \emph{Ondelettes et op\'erateurs, I: Ondelettes},
     Hermann, Paris, 1990.
      [English tanslation: \emph{Wavelets and operators}, Cambridge University Press
      1992.]
\bibitem{Pee}
\textsc{J. Peetre}, \emph{New Thoughts on Besov Spaces}, Duke Univ.
Math. Series, Durham, N.C., (1976).
\bibitem{Person}
\textsc{L. Person},
\emph{Interpolation with a parameter function}, Math. Scand. 59, (1986), 199--222.
\bibitem{Po}
 \textsc{W. Pompe},
    \emph{Unconditional biorthogonal bases in $L^p(\SR^d)$},
    Colloq. Math.,  92, (2002), 19--34.

\bibitem{Rou1}
\textsc{S. Roudenko}, \emph{Matrix-weighted Besov Spaces}, Trans.
Amer. Math. Soc., 355 (1), (2003), 273--314.

\bibitem{Rou2}
\textsc{S. Roudenko},
    \emph{The Theory of Function Spaces with Matrix Weights},
    Ph. D. Thesis, Departament of Mathematics, Michigan State
    University, (2002).


\bibitem{T1}
  \textsc{V.N. Temlyakov},
    \emph{The best m-term approximation and Greedy Algorithms},
    Advances in Comp. Math., 8, (1998), 249--265.

\bibitem{Wo1}
  \textsc{P. Wojtaszczyk},
    \emph{Wavelets as unconditional bases in $L\sb p(\SR)$},
    J. Fourier Anal. Appl. 5, no. 1, (1999),  73--85.
\bibitem{WoH}
\textsc{P. Wojtaszczyk},
\emph{Greediness of the Haar system in rearrangement invariant spaces},
in Approximation and Probability, Banach Center Publications, 72, T. Figiel and A. Kamont Eds., Warszawa (2006), 385--395.


\end{thebibliography}

\end{document}